\documentclass[12pt]{article}
\usepackage{fullpage,graphicx,psfrag,amsmath,amsfonts,verbatim}
\usepackage[small,bf]{caption}
\usepackage{amssymb}
\graphicspath{{figures/}}

\usepackage{listings}
\usepackage{longtable}
\usepackage{subcaption}
\usepackage{multirow}

\usepackage{color}

\definecolor{codegreen}{rgb}{0,0.6,0}
\definecolor{codegray}{rgb}{0.5,0.5,0.5}
\definecolor{codepurple}{rgb}{0.58,0,0.82}
\definecolor{backcolour}{rgb}{0.95,0.95,0.98}

\lstdefinestyle{mystyle}{ backgroundcolor=\color{backcolour},
      keywordstyle=\color{magenta}, numberstyle=\tiny\color{codegray},
      stringstyle=\color{codepurple}, basicstyle=\ttfamily\small,
      breakatwhitespace=false, breaklines=true, captionpos=t, keepspaces=true,
      numbers=left, numbersep=5pt, showspaces=false, showstringspaces=false,
      showtabs=false, tabsize=2, }

\definecolor{seagreen}{rgb}{0.18, 0.55, 0.34}
\definecolor{mediumviolet-red}{rgb}{0.78, 0.08, 0.52}
\definecolor{khaki}{rgb}{0.94, 0.9, 0.55}

\lstdefinelanguage{mypython}
{ keywords=[1]{from, import, as, assert, not, print, nonneg, PSD, axis},
      keywordstyle=[1]{\color{mediumviolet-red}}, keywords=[2]{cp, lo, pl,
      cvxpy, Variable, Parameter, sqrt, exp, numpy, np, Problem, Minimize,
      Maximize, value, solve, inner, sum, multiply, arange, range, norm1, norm2,
      norm_inf, abs, square, diagonal, outer, pos, hstack, power},
      keywordstyle=[2]{\color{seagreen}}, upquote=true, showstringspaces=false,
      basicstyle=\ttfamily, columns=fullflexible, keepspaces=true,
      emph={True,False,def,return,float,class,match,switch,len},
      emphstyle={\color{seagreen}}, belowskip=1em, aboveskip=1em,
      morecomment=[l]{\#} }

\newcommand{\BEAS}{\begin{eqnarray*}}
    \newcommand{\EEAS}{\end{eqnarray*}}
    \newcommand{\BEA}{\begin{eqnarray}}
    \newcommand{\EEA}{\end{eqnarray}}
    \newcommand{\BEQ}{\begin{equation}}
    \newcommand{\EEQ}{\end{equation}}
    \newcommand{\BIT}{\begin{itemize}}
    \newcommand{\EIT}{\end{itemize}}
    \newcommand{\BNUM}{\begin{enumerate}}
    \newcommand{\ENUM}{\end{enumerate}}

    \newcommand{\BA}{\begin{array}}
    \newcommand{\EA}{\end{array}}

    \newcommand{\eg}{{\it e.g.}}
    \newcommand{\ie}{{\it i.e.}}

    \newcommand{\reals}{{\mbox{\bf R}}}

    {\begin{quote}}{\end{quote}}

    \makeatletter
    \long\def\@makecaption#1#2{
       \vskip 9pt
       \begin{small}
       \setbox\@tempboxa\hbox{{\bf #1:} #2}
       \ifdim \wd\@tempboxa > 5.5in
            \begin{center}
            \begin{minipage}[t]{5.5in}
            \addtolength{\baselineskip}{-0.95pt}
            {\bf #1:} #2 \par
            \addtolength{\baselineskip}{0.95pt}
            \end{minipage}
            \end{center}
       \else
        \hbox to\hsize{\hfil\box\@tempboxa\hfil}
       \fi
       \end{small}\par
    }
    \makeatother

    \newcounter{oursection}

    \newcounter{lecture}

\usepackage{mathtools}
\usepackage{adjustbox}
\usepackage{booktabs}
\usepackage{csquotes}

\usepackage{hyperref}
\hypersetup{ colorlinks   = true,    
urlcolor     = blue,    
linkcolor    = blue,    
citecolor    = red      
}

\title{A Tax-Efficient Model Predictive Control Policy \\for Retirement Funding}
\author{Kasper Johansson \and Stephen Boyd }

\begin{document}
\maketitle

\begin{abstract}
The retirement funding problem addresses the question of 
how to manage a retiree's savings to provide her with a constant
post-tax inflation adjusted consumption throughout her lifetime.
This consists of choosing withdrawals and transfers from and between several 
accounts with different tax treatments,
taking into account basic rules such as required minimum 
distributions and limits on Roth conversions, 
additional income (such as Social Security payments),
liabilities, taxes, and the bequest when the retiree dies.

We develop a retirement funding policy in two steps.  In
the first step, we consider a simplified planning problem in
which various future quantities, such as the retiree's remaining 
lifetime, future investment returns, and 
future inflation, are known.  Using a simplified model of taxes, we pose
this planning problem as a convex optimization problem,
where we maximize the bequest subject to providing
a constant inflation adjusted consumption target.
Since this problem is convex, it can be solved quickly and reliably.

We leverage this planning method to form a 
retirement funding policy, \ie, an algorithm that determines
the actions to take each year, based on information known at that time.
Each year the retiree forms a new plan for the future years, using the
current account values and life expectancy, and optionally, 
updated information such as changes in tax rates or rules. 
The retiree then carries out the 
actions from the first year of the current plan.
This update-plan-act cycle is repeated each year, a general
policy called model predictive control (MPC).
The MPC retirement policy reacts to the effects of uncertain 
investment returns and inflation,
changes in the retiree's expected lifetime or external income 
and liabilities, and changes in tax rules and rates.

We demonstrate the effectiveness of the MPC retirement policy 
using Monte Carlo simulation, with 
realistic statistical models of investment returns and inflation.
By simulating thousands of trajectories of realized investment 
returns and inflation, and the retiree's lifetime,
we show that the MPC retirement policy works well,
delivering a higher bequest than the standard 4\% rule.
\end{abstract}

\clearpage
\tableofcontents
\clearpage

\section{Introduction}
Retirement planning involves managing uncertainties in investment returns,
inflation, longevity (\ie, lifetime of the retiree),
external income and liabilities, 
and possibly legislative changes to deliver financial well-being in later
life. At its core, the problem requires balancing income, contributions to
and withdrawals from various retirement accounts, 
and tax implications, under dynamic and evolving conditions. 
This paper addresses these challenges by formulating a
simple model predictive control (MPC) policy for retirement funding. 
The policy determines a retiree's withdrawals and transfers from and
between several accounts with different tax treatments, while taking into
account external income and liabilities, and the bequest when the retiree dies.
The policy is designed to deliver stable inflation adjusted consumption at a
target determined by the retiree, adapting the withdrawal and transfer plan
annually to changes in external factors such as life expectancy, realized investment
returns, and tax liabilities. Each year, the policy solves a convex optimization
problem to determine the optimal contributions, and withdrawals from the
accounts \cite{boydconvex}. The policy is simple, yet flexible, and can be
implemented using standard optimization software such as CVXPY \cite{cvxpy}. 
While we assume the target consumption is provided by the retiree, 
there are many ways to choose this target. 
For example, a retiree could use her current annual
consumption, or base it on the 4\% rule~~\cite{bengen1994determining}. She could
also use our method to choose the target consumption by simulating many
scenarios with different target values and choosing one that she is comfortable
with.
We provide an open-source implementation of the policy, along with all data needed
to replicate all numerical results given in this paper. 

Traditional approaches to retirement planning typically fall into two
broad categories, one mostly theoretical and the other mostly practical. Early
theoretical work such as Merton's framework for optimal
consumption~\cite{merton1975optimum}, offers elegant solutions based on
simplified assumptions. The Merton problem, a well-known stochastic control
model, focuses on determining the optimal consumption and investment strategy to
maximize expected constant relative risk aversion (CRRA) utility over time.
While it addresses uncertainty in portfolio returns, its original formulation
overlooks practical
considerations like longevity risk, taxes, the specific rules governing
different retirement accounts, and the fact that a typical objective of a
retiree is stable inflation adjusted 
consumption throughout retirement, rather than maximizing an expected utility.
Recent extensions of the Merton problem have
considered uncertainty in longevity~\cite{moehle2021certainty}. 
In Merton's problem, and its extensions, consumption is handled as a 
variable, and judged by an increasing concave utility function.  In contrast,
we strive for inflation adjusted consumption that is at, or near,
a constant specified target value.

At the other extreme, financial advisors use detailed spreadsheets
and customized heuristics to guide clients. They assess a client's financial
situation, future income needs, and risk tolerance to craft personalized
strategies, which typically involve selecting appropriate retirement accounts
(\eg, 401(k)s, IRAs, annuities), diversifying investments, and offering advice
on tax-efficient withdrawals, Social Security benefits, and healthcare costs.
They also monitor portfolios and adjust strategies in response to market shifts,
life events, and changing goals.

We propose a third approach: an MPC-based framework that
solves an optimization problem each year to create a plan for future actions,
based on current projections for investment returns and inflation, 
retiree lifetime, and other uncertain future quantities.  This method incorporates
the retiree's age, financial situation, income streams, liabilities, and
projections for longevity, investment returns, and inflation.
Our model explicitly accounts
for the main tax rules as well as contribution and withdrawal restrictions associated
with different account types. By recalibrating the plan annually based on
updated information, our framework adapts to changes in 
both forecasted uncertainties and
realized outcomes, offering a dynamic and robust strategy for retirement
funding. Notably, the effectiveness of our approach does not hinge heavily on
the precision of the forecasts or the specific formulation of the planning
problem, as the iterative nature of the framework inherently adjusts to evolving
conditions and uncertainties~\cite{mcallister2022stochastic, grimm2007nominally,
kuntz2024beyond}. In essence, we propose an algorithmic method 
to carry out retirement funding.
Our method models only the main tax rules, and makes a number of 
approximations; it could be monitored or
modified by a financial retirement advisor who can make adjustments
depending on the specific needs and situation of a client.

\subsection{Related work}
We mention here some related work, ranging from portfolio
construction and tax-aware investment strategies to pure retirement planning.

\paragraph{Tax-aware investing.}
Our approach is related to tax-aware investment
strategies~\cite{constantinides1983capital, constantinides1984optimal,
moehle2021tax}. These papers consider portfolio construction and asset
allocation, taking taxes into account. The underlying
idea is that investors can reduce their tax liability by deferring capital gains
until they are classified as long term, and then taxed at a lower
rate, and realizing short term losses. (While not addressed in these papers,
bequests of assets with embedded gains also have beneficial tax treatments.)
While the simplest version of our retirement policy does
not consider portfolio construction or optimal tax-aware rebalancing, it could
certainly be combined with such strategies.

\paragraph{Portfolio construction.} Portfolio construction refers to the process
of selecting a combination of assets that optimally balances the trade-off
between expected return and risk.
The foundation of modern portfolio construction dates back 
to Markowitz's seminal 1952
paper~\cite{markowitz1952portfolio}, introducing the concept of modern portfolio
theory, which models portfolio returns as random variables. The key innovation
was the formulation of the trade-off between expected return and portfolio risk,
with the goal of maximizing the portfolio's expected return for a given level of
risk. 
Numerous advancements have been built on Markowitz's initial 
formulation~\cite[Chap.~14]{grinold2000active},
\cite[Chap.~6]{narang2024inside}, \cite{cornuejols2006optimization, kolm201460,
chakrabarty2023mathematical, gunjan2023brief, boyd2024markowitz}. 
Convex optimization, which we use in our retirement policy, also plays
a large role, and is widely used, in portfolio construction.

This paper does not focus on portfolio construction per se. Instead we fix
the investment portfolios for the different accounts to be reasonable simple ones.
The method we describe could certainly be combined with more sophisticated
portfolio construction methods, such as ones that include collar options to 
limit investment returns.
The choice of portfolios is reflected in our retirement plans in a simple
way, as a fixed or time varying typical return for each account.
When we assess our retirement policy, however, we use realistic statistical models
of investment returns, inflation, and longevity.

\paragraph{Retirement funding.} The key challenge in retirement funding is
determining how to invest and spend savings to maintain a comfortable lifestyle
throughout the retiree's later years. The literature on retirement funding is
vast, with numerous approaches and strategies. A simple strategy is proposed
in~\cite{bengen1994determining}, where the author argues, based on historical
data, that a retiree can safely withdraw 4\% of their initial savings each year,
adjusted for inflation. Much research has focused on these fixed withdrawal
strategies (with adjustments for inflation), and the objective is often to find
a strategy with low probability of running out of money before the end of
life~\cite{blanchett2023redefining}. Other studies advocate variable consumption
in proportion to survival probabilities, adjusted up for exogenous pension
income and down for longevity risk aversion~\cite{milevsky2011spending}, or
allow the withdrawal to vary, depending on the retiree's sensitivity to
consumption changes, arguing that wealthy individuals are less sensitive to
spending cuts~\cite{idzorek2019ldi, blanchett2023redefining}. Our simple MPC
policy has elements of both of these approaches. Given a target consumption, it
adjusts the withdrawals and transfers between accounts to give an inflation
adjusted consumption that with very high probability is close to the target,
and a bequest that is as large as possible.


\paragraph{Model predictive control.}
Our method is based on a well known idea, used in many other
fields. The goal is to make sequential decisions (that affect the future), 
despite uncertainty~\cite{kochenderfer2015decision}.
In MPC this is done in a simple way that at each decision time carries
out three steps. In the first step the uncertain future quantities are
predicted or forecasted.  In the second step, an optimal sequence of actions,
or action plan, is found using an optimization algorithm, 
assuming that the forecasts are correct.
In the third step, we take the action that is the first action in the
plan.  The whole forecast-plan-act process repeats at the 
next decision time, 
with updated forecasts, and taking into account the current state.
For general background on MPC, see, \eg, 
\cite{borrelli2017mpc} and the references therein.

MPC is used in a variety of industries, sometimes under different names such as
receding horizon control, dynamic matrix control, or shrinking horizon control.
Its applications are broad, and include chemical
engineering~\cite{kumar2012model}, supply chain management~\cite{wang2007model,
perea2003model}, multi-period trading~\cite{li2022multi, BoydKahnMultiPeriod},
optimal execution of an order~\cite{boyd2024cvx}, and even the autonomous
landing of SpaceX's Falcon 9 and Falcon Heavy first 
stages~\cite{blackmore2016autonomous}. 
These applications often rely on convex
optimization, which ensures that the resulting problems are tractable and
can be efficiently and reliably solved \cite{boydconvex}.

\subsection{Contributions}
We make two main contributions. First, we formulate the retirement funding planning
problem as a mathematical optimization problem, which is convex and so 
can be solved reliably and efficiently.  
(This also means that we can evaluate the policy over many
thousands of scenarios, which is necessary to assess its effectiveness 
by simulation.)
Second, we leverage the planning problem into a
simple MPC policy, which systematically adapts the
retiree's retirement plan to changing external factors such as life expectancy,
investment returns, and tax liabilities. 

\subsection{Outline}
In~\S\ref{sec:retirement} we outline the retirement planning problem,
making a number of strong assumptions such as knowing future investment 
returns, inflation, and the retiree's lifetime.
The method does, however, take into account the main rules
that govern distributions and transfers between accounts, as well as a basic
model of taxation.
We formulate this simplified planning problem as a convex optimization
problem, which can be very efficiently solved.

In \S\ref{sec:mpc} we leverage this planning method into a retirement 
funding policy,
which each year suggests deposits to,
withdrawals from, and transfers between accounts, based on 
information available at that time.
Because we re-plan each year, this policy reacts to uncertainties in
investment returns, inflation, retiree lifetime, external income and liabilities,
and even changes to the rules and tax rates.

In \S\ref{sec:data_and_models} and \S\ref{sec:experiments} we
evaluate our proposed retirement funding policy.
We introduce statistical models for investment returns, Treasury rates, and
inflation in \S\ref{sec:data_and_models}, describing both how to use 
these to simulate realistic scenarios and to form forecasts for MPC. 
We evaluate the MPC policy using simulations 
in~\S\ref{sec:experiments}.
Extensions and variations are discussed in~\S\ref{sec:extensions}, and we
conclude in~\S\ref{sec:conclusion}.
In appendix~\ref{sec:collar_investments} we give numerical results
for our method using more sophisticated portfolios than a simple stock/bond 
split.

\section{Optimal retirement planning}\label{sec:retirement} 
In this section we introduce our notation and
describe the retirement planning problem.
We consider a US tax resident who retires in year $t=1$, 
and uses her savings (and any additional
income) to fund her retirement consumption for years $t=1, \ldots, T$.
We assume the retiree retires at age 60 or older, which simplifies some
of the rules and eliminates some penalties.
All USD values are in real dollars, \ie, inflation adjusted.

In year $T$ the retiree dies, and leaves a bequest.
Assuming we know the year of death of the retiree is clearly 
unrealistic; this will be addressed in \S\ref{sec:mpc}, where 
we use the retirement plan developed in this section to create a 
policy.

\subsection{Savings accounts}
The retiree has three types of accounts, which are taxed differently: a
brokerage account, for which realized gains are taxed at the long-term capital
gains rate (\ie, assuming they have been held for at least one year),
a traditional IRA, for which withdrawals are taxed, and a Roth IRA, where
withdrawals are tax-free. (We do not consider additional accounts, such as a
401(k) or Roth 401(k) here, but discuss how to model them
in~\S\ref{sec:extensions}.) 

We denote the (nonnegative) amounts in these
accounts at the beginning of year $t$ as $B_t$, $I_t$, and $R_t$, respectively,
for $t=1, \ldots, T+1$, denoted in real USD, \ie, adjusted for inflation. We
assume the retiree starts off with balances $B_1=B^{\text{init}}$,
$I_1=I^{\text{init}}$, and $R_1=R^{\text{init}}$. To finance retirement in year
$t=1,\ldots,T$ the retiree withdraws amounts $b_t$, $i_t$, and $r_t$, from the
brokerage, IRA, and Roth accounts, respectively, with $b_t \leq B_t$, $i_t \leq
I_t$, and $r_t\leq R_t$. These can be negative, which corresponds to a deposit.

\subsection{Distributions}
This section describes the rules for distributions
from brokerage, IRA, and Roth
accounts, including components like conversions, deposits, and withdrawals. It
covers withdrawal limits, deposit restrictions, and required minimum
distributions (RMDs) for IRAs, as well as conditions for Roth conversions.

\paragraph{Brokerage withdrawals.} Each year the retiree withdraws an amount 
$b_t$ from the brokerage account, with $b_t \leq B_t$. 
This amount can be negative, which
means a deposit of $-b_t$ to the brokerage account. 

\paragraph{IRA withdrawals.} The withdrawal from the IRA account is split
into
three components, $i_t^c$, $i_t^d$, and $i_t^w$, all nonnegative. 
The first component $i_t^c$ is a conversion to the Roth
account, the second component $i_t^d$
 is the amount we deposit into the IRA account, and
the third component $i_t^w$ is the withdrawal. The net IRA withdrawal is
\[
i_t = i_t^c - i_t^d + i_t^w.
\]
This number can be negative, indicating a net deposit into the IRA; however,
this is unlikely for most retirees following the RMDs,
described below. We split the net withdrawal into these three components to
handle Roth conversion limits and RMDs.


\paragraph{Roth withdrawals.} 
The Roth withdrawal is also
split into three components, $r_t^c$, $r_t^d$, and $r_t^w$,
all nonnegative. The first component, $r_t^c$, is a conversion from the IRA
account, the second, $r_t^d$, is the amount we deposit in the Roth account, and
the third is the withdrawal. The net Roth withdrawal is
\[
r_t = -r_t^c - r_t^d + r_t^w.
\]
The net Roth withdrawal can be negative, indicating a net deposit into the
Roth account, most commonly attained through a Roth conversion. 

There are some limits on Roth withdrawals we do not account for here.
Roth withdrawal limits mainly affect individuals under the age of 59.5, and so
will most likely not be relevant for a retiree. We discuss this further in
\S\ref{sec:extensions}.

\paragraph{Deposit limits.}  
The total annual deposits to traditional IRA and Roth IRA accounts are subject
to certain restrictions:  
\[
i_t^d + r_t^d \leq d^{\max},
\]  
where $d^{\max}$ depends on factors such as age and income. For 2024,
individuals aged 50 or older with a modified adjusted gross income (MAGI) below
\$146,000 are allowed to contribute up to \$8,000. (Roth IRA contributions have
income limits, but deductible traditional IRA contributions do not, unless the
individual is covered by a workplace retirement plan.) A more detailed
description of MAGI dependences is given in~\S\ref{sec:extensions}. 

Additionally, contributions to an IRA or Roth IRA must come from earned income
during the year, expressed as
\[
i_t^d + r_t^d \leq e_t,
\]  
where $e_t$ represents earned income in year $t$. Combining these constraints,
the total contributions are bounded by
\[
i_t^d + r_t^d \leq \min\{d^{\max}, e_t\}.
\]

\paragraph{Required minimum distributions and Roth conversions.}
The IRA account is subject to RMDs. This states
that $i^w_t \geq \kappa_t I_t$, where $\kappa_t \in [0,1]$ is the given fraction
of the account value that must be withdrawn in year $t$. These fractions are
zero until the retiree reaches age 73, and after that given as $1/y_t$, where
$y_t$ is the expected remaining lifetime of someone alive in year $t$,
obtained from tables published by the Internal Revenue Service (IRS).

Roth conversions are also limited. When $r^c_t>0$, the amount $r^c_t$ must come
from an IRA withdrawal in excess of the RMD. We can express the RMD and Roth
conversion constraints as
\[
i^w_t \geq \kappa_t I_t, \quad i^c_t = r^c_t, \quad t=1, \ldots, T.
\]
This states that the IRA withdrawal must be at least the RMD, and that the Roth
conversion must come from an IRA withdrawal in excess of the RMD, \ie, from the
portion of the net IRA withdrawal that corresponds to the conversion.


\subsection{Income, liabilities, and consumption}
We now outline the retiree's sources of income and liabilities, including earned
income, such as a salary, which can be used for funding savings accounts or
consumption, as well as additional income like Social Security or pension
payments. The section also addresses liabilities, which represent the retiree's
known future financial obligations.

\paragraph{Earned income.}
If the retiree has earned income, such as a salary, she may use this to fund any
of her accounts, or her consumption. We denote the (nonnegative) 
amount of earned income in year $t$ as $e_t$.

\paragraph{Additional income.}
The retiree may receive additional income, such as Social Security payments,
pension income, or annuity payments, which falls outside the category of earned
income. We denote this income as $a_t$, representing the amount received in year
$t$. Additional income cannot be used to fund the IRA or Roth accounts.

\paragraph{Liabilities.}
Liabilities include, \eg, mortgage payments or other debts. We denote the
liability in year $t$ as $l_t$. The liability is typically nonnegative, but
could be
negative, meaning an income $l_t$. We distinguish between a negative liability
and additional income, as we pay taxes on the latter, but not on the former. 
An example of a negative liability is a tax refund, which we do not pay
taxes on.

\paragraph{Constant consumption target.}
We assume that the retiree wants to consume a constant nonnegative amount 
$c$ in each year. Since we work with real (inflation adjusted) values, this
corresponds to consumption that includes a cost of living adjustment (COLA) 
each year. The goal is for this consumption to be close to a target value,
which we denote as $c^{\text{tar}}$. 
The target level $c^{\text{tar}}$ can be chosen in various ways, such as
the retiree's average consumption over recent (pre-retirement) years, 
or applying a rule of thumb like the 4\% rule.

\paragraph{Taxes.}
The retiree pays a nonnegative tax bill in year
$t=1,\ldots, T$, based on taxable income and capital gains 
in year $t$, denoted $\tau_t$. 
(This amount is described in more detail below.)

\paragraph{Cash balance.}
Cash balance must hold, which means
\[
b_t + i_t + r_t + e_t + a_t = c + l_t + \tau_t, \quad t=1, \ldots, T.
\]
The lefthand side is the net amount withdrawn from the three accounts plus
income; the righthand side is the retiree's total expenses, including
consumption, taxes, and liabilities. 

\subsection{Taxes}

\paragraph{Taxable income.}
The taxable income in year $t$ is given by
\[
\omega_t = i^c_t - i^d_t + i^w_t + e_t + a_t.
\]
The first three terms correspond to the net IRA withdrawal; we pay taxes on the
Roth conversion and the withdrawal, while the deposit acts as a tax deduction.
The last
two terms are the earned income and additional income. 

\paragraph{Capital gains.}
Capital gains in the brokerage account are taxed
separately from income and are based on nominal gains, the difference
between the sale and purchase price in nominal (not inflation adjusted) dollars.
(With this one exception, all other dollar values in this paper are 
inflation adjusted.)
To track taxable gains in the brokerage account, we record tax lots, maintaining
the number of shares $s_t$ and nominal purchase price $p_t$ at the time of
purchase. We assume shares are sold using the average cost basis method, where
the cost basis per share is given by  
\[
x_t = \frac{\sum_{\tau < t} s_\tau p_\tau}{\sum_{\tau < t} s_\tau},
\]  
denoted in nominal USD per share.
Assuming no tax benefits from losses, the capital gain (in real dollars) upon sale at 
time $t$ is  
\[
\omega^\text{cg}_t = (b_t)_+ \left(1 - \delta_t\right), 
\]  
where $(u)_+ = \max\{u,0\}$ and
$\delta_t = x_t/p_t$ is the ratio of basis to value.
Since future sales prices and cost basis are unknown when constructing the plan,
we approximate this ratio using the last known value $\delta_0 =
x_0/p_0$, leading to the simplified expression  
\[
\omega^\text{cg}_t = (b_t)_+ \left(1 - \delta_0\right), \quad t = 1, \dots, T.
\] 
Thus our approximation is that the ratio of basis to value in the 
brokerage account remains constant.

We assume positive capital gains are taxed at a fixed nonnegative rate $\xi$, 
reflecting the retiree's long-term capital gains rate, so the 
capital gains tax (in real dollars) in year $t$ is given by $\xi \left(\omega^\text{cg}\right)_+$.
Here we ignore the tax reduction possible when a withdrawal from the brokerage account
results in a capital loss, which occurs when $\delta_0 > 1$.
 
\paragraph{Tax bill.}
The tax bill is given by
\[
\tau_t = \phi(\omega_t) + \xi (b_t)_+ \left(1 - \delta_0\right)_+,
\]
where $\phi$ is an increasing convex piecewise linear function with
$\phi(\omega)=0$ for $\omega \leq 0$. (For simplicity, we assume there is no tax
benefit from a negative taxable income.) The tax function is described by the
$K$ kink points $0<\beta_1 < \cdots < \beta_K$, and the slopes of $\phi$ in the
$K+1$ intervals $[0,\beta_1], [\beta_1,\beta_2],\ldots, [\beta_{K-1},\beta_K],
[\beta_K,\infty)$, given by $0< \eta_1 < \cdots <\eta_{K+1}$. The numbers
$\beta_i$ are the tax bracket boundaries, and the slopes $\eta_i$ are the
marginal tax rates in the $K+1$ tax brackets.
If we include the kink in $\phi$ at $0$, we have $K+1$ kink points, and $K+2$
tax brackets, with zero tax for the negative income tax bracket $(-\infty,0]$.

We model the tax function as constant in future years. Since we work with 
inflation adjusted values, this corresponds to annual adjustments of the 
tax brackets in future years, which is a reasonable assumption.

\paragraph{Tax simplifications.}
We have made several simplifying assumptions. First, we assume that the retiree
pays a fixed tax rate $\xi$ on capital gains, rather than distinguishing between
short-term and long-term gains and applying a progressive tax rate based on
taxable income. Moreover, while long-term capital gains are taxed at a lower
rate than ordinary income at the federal level, most states treat capital gains
as ordinary income. Second, we assume that the retiree receives no tax benefits
from capital losses or negative income. 
These assumptions simplify the problem, and introduce some
conservatism in the planning. Third, we have ignored dividends, which are taxed
annually, either as ordinary income or at the long-term capital gains rate,
depending on whether the dividend is qualified and the holding period.  
Fourth, we assume that the
retiree pays full taxes on the additional income, which is often not the case.
For example, Social Security payments are only partially taxed, somewhere
between 0\% and 85\%, depending on the retiree's income. Our simplification that
additional income is fully taxed also introduces conservatism. We can modify the
problem to account more accurately for these effects;
see~\S\ref{sec:extensions}. 

These simplifications are good enough for the purposes of creating a funding
policy with MPC.
When we \emph{simulate} the retirement funding policy, however,
we compute the capital gain tax bill exactly
using the progressive tax rate, which depends on the retiree's taxable income.

\subsection{Account dynamics}
This section outlines the dynamics of the investment accounts, including 
investment gains and a bequest. 
Account balances evolve based on inflation adjusted returns, and
the bequest is defined as the sum of the account values at the end of the
investment horizon, \ie, any remaining funds after the retiree's death.

\paragraph{Investment gains.}
The account balance dynamics are given by
\[
B_{t+1} = (B_t-b_t) \rho^B_t, \quad
I_{t+1} = (I_t-i_t) \rho^I_t, \quad
R_{t+1} = (R_t-r_t) \rho^R_t, \quad t=1, \ldots, T,
\]
where $\rho^B_t, \rho^I_t, \rho^R_t$ are the (positive, inflation adjusted) investment
total returns for the three accounts. The inflation adjusted returns are in
general not known, but would be estimated based on the investment portfolio and
current market conditions. For example, if the retiree holds a risk-free bond
portfolio paying 4\% annually, and the inflation is 2\%, then the inflation
adjusted return would be 2\%.

The assumption of known investment returns is of course unrealistic.
But it is adequate for the purpose of making a plan to be used in an MPC
policy.  Here too we can introduce conservatism in the planning process,
by taking values of the constant returns $\rho^B, \rho^I, \rho^R$ that 
are lower than we might really expect.  For example, they might be taken
as the 30th or 40th percentile of anticipated (random) investment returns,
as opposed to a mean or median. When we carry out simulation, however,
we update the account balances
using realized returns, as described in~\S\ref{sec:mpc}.

\paragraph{Bequest.}
The retiree dies in year $T$ and leaves the bequest
\[
q = B_{T+1} + I_{T+1} + R_{T+1}
\]
in year $T+1$.
Knowing the year of death of a retiree is also an unrealistic assumption, 
like knowing future investment returns or inflation,
but it is also adequate for the purposes of MPC.  
Like the investment returns, we can introduce conservatism in the plan
by over-estimating the remaining lifetime of the retiree, \eg,
using the 75th percentile of remaining lifetime from a life table, or 
increasing expected lifetime by a fixed factor. 

\subsection{Objective}
The primary objective of the retiree is to maximize the bequest, subject to
maintaining a constant consumption close to the target value, if possible.
We model this objective as
\[
U(c,q) = q - \gamma (c^{\text{tar}}-c)_+,
\]
where $\gamma > 0$ represents the retiree's preference for maintaining a constant
consumption close to the target value. With a suitable (large enough) choice of
$\gamma$, this objective incentivizes choosing
$c=c^{\text{tar}}$ if at all possible.
The parameter $\gamma$ can be chosen, for example, by simulating the retirement
funding policy with multiple choices, and choosing one that gives good
performance. 
The second term, which is proportional to the consumption shortfall
$(c^\text{tar}-c)_+$, can be thought of as a softening of the hard constraint 
$c=c^\text{tar}$ \cite{boyd2024markowitz}.  Roughly speaking, we hope that 
in all but very adverse circumstances, we have $c=c^\text{tar}$.
That is, we reduce consumption below the target only if we really must.



\subsection{Retirement funding planning problem}\label{e-rfp}
We formulate the choice of an optimal plan for retirement funding
as an optimization problem, with the goal of
maximizing our objective subject to the many constraints described above.
This results in the optimization problem
\BEQ\label{prob:nonconvex}
\begin{array}{ll}
\mbox{maximize} & U(c, q) \\
\mbox{subject to} & B_1 = B^{\text{init}}, \quad I_1 = I^{\text{init}}, \quad
R_1 = R^{\text{init}}, \\ 
& B_t \geq 0, \quad I_t \geq 0, \quad R_t \geq 0, \quad t = 1, \ldots, T+1,\\
& B_{t+1} = (B_t - b_t) \rho^B_t, \quad I_{t+1} = (I_t - i_t) \rho^I_t, \quad
R_{t+1} = (R_t - r_t) \rho^R_t, \quad t = 1, \ldots, T,\\
& i_t = i_t^c - i_t^d + i_t^w, \quad r_t = -r_t^c - r_t^d + r_t^w, \quad t=1,
\ldots, T,\\
& b_t \leq B_t, \quad i_t \leq I_t, \quad r_t \leq R_t, \quad t = 1, \ldots,
T,\\
& b_t + i_t + r_t + e_t + a_t = c + l_t + \tau_t, \quad t=1, \ldots, T, \\
& i^w_t \geq \kappa_t I_t, \quad i^d_t + r^d_t \leq \min\{d^{\max}, e_t\},  \quad
i^c_t = r^c_t, \quad t=1, \ldots, T,\\
& \tau_t = \phi(i^c_t - i^d_t + i^w_t + e_t + a_t) + \xi (b_t)_+ (1 - \delta_0)_+,
 \quad t = 1, \ldots, T, \\
& q = B_{T+1} + I_{T+1} + R_{T+1}.
\end{array}
\EEQ 
We refer to this as the retirement funding planning (RFP) problem.
The variables and parameters are summarized below.

\paragraph{Variables.}
The variables are \BIT
\item \emph{Account balances} $B_t, I_t, R_t$, $t=1, \ldots, T+1$.
\item \emph{Distributions} $b_t, i_t, r_t, i^c_t, i^d_t, i^w_t, r^c_t, r^d_t,
r^w_t$, $t=1,\ldots, T$.
\item \emph{Consumption and bequest} $c, q$.
\item \emph{Tax bill} $\tau_t$, $t=1, \ldots, T$. 
\EIT

\paragraph{Data.}
The given data are \BIT
\item \emph{Initial account balances} $B^{\text{init}}, I^{\text{init}},
R^{\text{init}}$.
\item \emph{Ratio of brokerage basis to value} $\delta_0$.
\item \emph{Investment returns} $\rho^B_t,\rho^I_t,\rho^R_t$, $t=1, \ldots, T$.
\item \emph{Income and liabilities} $e_t, a_t, l_t$, $t=1, \ldots, T$.
\item \emph{RMD rates} $\kappa_t$, $t=1, \ldots, T$.
\item \emph{Tax brackets and marginal tax rates} $\beta_1, \ldots, \beta_K$,
$\eta_1,\ldots, \eta_{K+1}$.
\item \emph{Long-term capital gain tax rate} $\xi$.
\item \emph{Deposit limit} $d^{\max}$.
\item \emph{Target consumption} $c^\text{tar}$.
\item \emph{Objective trade-off parameter} $\gamma$. 
\EIT

\paragraph{Convexity.}
As stated, the RFP problem \eqref{prob:nonconvex} is not a convex 
optimization problem, due to the non-affine tax function $\phi$. 
However, it can be transformed to an equivalent convex optimization problem by
replacing this nonconvex constraint with its relaxation
\BEQ\label{e-relax-tax}
\tau_t \geq \phi(i^c_t - i^d_t + i^w_t + e_t + a_t) + \xi (b_t)_+ (1 - \delta_0), \quad t = 1, \ldots, T.
\EEQ
This results in a convex optimization problem~\cite{boydconvex}.

A standard argument shows that if we solve the RFP problem \eqref{prob:nonconvex}
with the relaxed tax constraints \eqref{e-relax-tax} (with inequality constraints),
replacing the orignal tax constraints (with equality constraints),
we obtain the solution of the original
problem, \ie, the tax constraints hold with equality.
It is easy to see why the relaxation is tight at the solution.
The relaxation means that we are given the option to pay more in taxes than what we
are required to, which is not something we would do; we could always reduce the tax
payment and use the additional money to increase $q$. Hence, the solution to the
convex problem is the same as that of the original nonconvex
problem~\eqref{prob:nonconvex}.

\paragraph{Solve speed.}
We henceforth refer to the relaxed RFP problem, which is \eqref{prob:nonconvex}
with the tax constraints relaxed as in \eqref{e-relax-tax},
as the RFP problem. It has around $13 T$
variables, $7 T$ equality constraints, and $9 T$ inequality constraints. 
It can be solved very reliably and quickly. 

As an example, consider a planning horizon of $T=45$
(\eg, age range 60--105), which yields an RFP problem with around $600$
variables, $300$ equality constraints, and $400$ inequality constraints. This
problem can be solved in around $0.01$ seconds on an Apple M2 Max chip (12-core
CPU, 64GB unified memory) using CVXPY and Clarabel~\cite{cvxpy,
goulart2024clarabel}. The solver is single-threaded, so we can solve around
500 tax planning problems per second, using all 12 cores. 
Our ability to solve the tax planning problem this quickly is evidently
not needed to implement MPC, where it will be solved once, or a few times, 
per year.  But the speed is critical 
when we carry out Monte Carlo simulations to evaluate the 
effectiveness of a policy.

\paragraph{Simplifying assumptions.}
We have already mentioned some strong simplifying assumptions made
in the formulation of the RFP, such as knowing the remaining lifetime
of the retiree, knowing future investment returns in 
her savings accounts, assuming gains in the brokerage account
are taxed at a fixed rate, and assuming that all external income
is taxed.

Other simplifications regard the IRA and Roth accounts. We assume that the
retiree can withdraw any amount from the IRA and Roth accounts, which is almost
accurate for those over 59.5 years old. In reality, withdrawals from Roth IRAs
are subject to the five-year rule. This rule requires that the account be open
for at least five years for earnings to be withdrawn tax-free, even after age
59.5. Furthermore, for Roth conversions, the five-year rule applies separately
to each conversion, requiring that five years have passed since the conversion
for the earnings to be withdrawn penalty-free. (There are
also penalties for withdrawals from IRAs before the age of 59.5, which we have
not considered.) These conditions could influence the retiree's withdrawal
strategy, particularly for those who rely heavily on Roth funds. Furthermore,
many retirees have additional accounts, such as a 401(k) or Roth 401(k), which
have different rules than IRAs, mainly higher contribution limits and possible
company matching. 

Our RFP problem can be extended to handle many of these aspects,
some discussed in \S\ref{sec:extensions}.
Our simulation experiments suggest that this level of accuracy in the
planning problem is not needed, at least 
when the RFP is used in an MPC policy.

\section{Model predictive control policy} \label{sec:mpc}

MPC is a dynamic optimization approach that determines a new plan of action each
year based on the current account balances (which depend on realized investment
returns) and updated forecasts, such as the retiree's expected year of death,
inflation rates, and investment returns. Essentially, the retiree solves the
planning problem~\eqref{prob:nonconvex} annually, with several modifications to
account for evolving conditions and changing values. We describe these in a bit
more detail below.

\subsection{Expected remaining lifetime}
The planning horizon $T$ is adjusted each year to reflect the retiree's
remaining expected lifetime. The current year $t=1$ always corresponds to the
present year, while the horizon $T$ is updated using life expectancy tables
such as those
provided by the Social Security Administration~\cite{ssa_actuarial_table}. For
example, a 65-year-old female retiree initially has a planning horizon of
approximately 20 years. This planning horizon decreases annually as 
the retiree ages.

We can modify this simple method to introduce conservatism in the planning
process.  Instead of adjusting $T$ to be the expected remaining lifetime, we can
take it to be a quantile such as the 70th percentile, according to the
distribution of lifetimes, or increase it over the expected lifetime 
by some fixed factor such as 30\%.
Roughly speaking, we plan on the (optimistic)
basis that the retiree's remaining lifetime will be longer than the mean.
This introduces some conservatism, since we must plan on funding
the retiree's consumption over a longer period. 

\subsection{Updated account values and forecasts} \label{sec:updates}
The account values $B^{\text{init}}, I^{\text{init}},$ and
$R^{\text{init}}$ are set to be the current values of the accounts.
These updates incorporate realized investment returns, taxes,
deposits, withdrawals, external income, and liabilities, over the last year.

We can also, optionally, update our forecasts and even models.
For example, the tax rates and brackets in our planning problem might 
be the current ones, implicitly indexed for inflation since we work with 
inflation adjusted quantities.  These could be replaced by the 
updated true current values, or an entirely new model of tax rates,
if they are changed by new tax laws.
Other examples include forecasts of inflation or investment returns, 
earned or additional income, and liabilities.
These updates modify the problem data to reflect the current conditions.

\subsection{Cash balance discrepancy} \label{sec:cash_discrepancy}  
The cash balance equation used in the MPC policy,  
\[
b_t + i_t + r_t + e_t + a_t = c + l_t + \tau_t, \quad t=1, \ldots, T,
\]
relies on forecasts of uncertain quantities. For instance, income streams and
liabilities may not be precisely known at the beginning of the year, and the tax
bill $\tau_t$ depends on realized taxable income.  In practice, capital gains
are taxed progressively, with tax rates ranging from 0\% to 20\% depending on
the retiree's total taxable income. However, in the planning problem, we assume
a fixed tax rate $\xi$, which can lead to discrepancies between the planned
and actual cash balance.  

To account for this discrepancy, we introduce a realized adjustment term,  
\[
\Delta_t = c + \tau_t + l_t - (b_t + i_t + r_t + e_t + a_t),
\]
which represents the difference between actual and planned cash flows. This
discrepancy is carried over to the next year's liabilities by adjusting
$l_{t+1}$. For example, if the retiree underestimates her tax bill, the
shortfall is treated as an additional liability in the following year.

\subsection{Simulation}
When simulating a policy, we need to account for the fact that some parameters
in the planning RFP problem are not known exactly, as discussed in
\S\ref{sec:updates} and \S\ref{sec:cash_discrepancy}.
To account for this (generally small) discrepancy, we carry out the simulations
as follows.
At the beginning of each year we solve the RFP problem in MPC to determine the 
current year's deposits and withdrawals from the savings accounts.
This is done using forecasted investment returns.
At the end of the year we update the account values using the realized investment 
returns, rather than the forecasted returns.
The realized investment
returns are used to update the account values, and the tax
bill is computed based on the realized taxable income, and the capital gains
tax, which can differ from the forecasted tax bill 
that appears in the RFP solved at the beginning of the year.
Next, we perform cash
accounting by adding the cash balance discrepancy to the next year's
liabilities, as described in~\S\ref{sec:cash_discrepancy}. (When the retiree
dies, we subtract any remaining liabilities from the bequest.)

\subsection{Evaluating the MPC policy}
We start by solving a big problem (as above), \ie, one with large $T$;
the problem gets smaller (\ie, $T$ shrinks) every year until the retiree dies.
For example, consider a 60-year-old retiree with an original
45-year planning horizon. 
The first plan is the problem example described earlier,
which can be solved in around $0.01$ seconds using CVXPY and 
Clarabel~\cite{cvxpy, goulart2024clarabel}.
Subsequent planning problems are smaller (since $T$ decreases from $45$),
and can be solved faster.
Assuming the retiree lives to be 105, we solve 45 planning problems,
requiring a total time of around 0.13 second, or about 0.003 seconds per
problem, with Clarabel.
The retiree will likely have a shorter lifespan, which means
that evaluating the MPC policy for one retiree's lifetime 
can be simulated in even less time. This means that by multi-threading, say
using all 12 cores of the Apple M2 Max chip, we can evaluate roughly 5,000 
retirement policies in under a minute. This only takes into account the
time to solve the optimization problems. The actual simulation can be up to 
five times longer due to our inefficent implementation of the Python code 
for simulation.  An efficient implementation would have a total simulation
time closer to the solve time, around one minute.

\section{Data and forecasts} \label{sec:data_and_models}

In this section we describe our statistical model of inflation adjusted 
investment returns, Treasury rates, and inflation,
which we use to evaluate the MPC policy.
We construct simple statistical models for investment returns, Treasury rates, 
and inflation rates.  We use these models
to generate realistic data for use in Monte Carlo simulations,
and also to create simple forecasts of future values for use in MPC.

We focus here on traditional stock/bond investment portfolios since 
they are simple and commonly used.  More sophisticated portfolios
that limit the left and right tails of the return distribution using collars, 
Treasury inflation securities (TIPS), and other methods 
are discussed in appendix~\ref{sec:collar_investments}.

\subsection{Market returns}  
\begin{figure}
\centering
\includegraphics[width=0.6\textwidth]{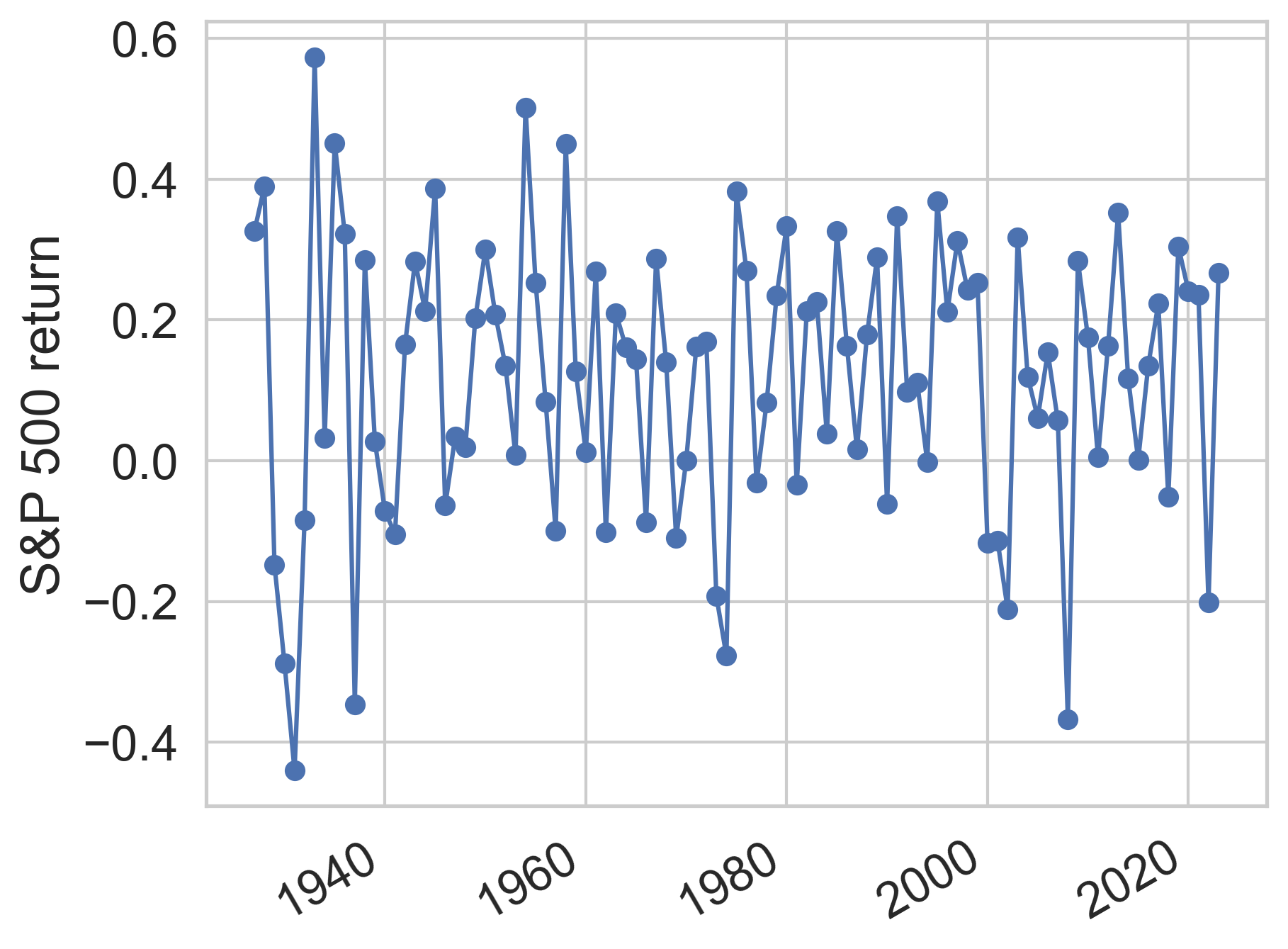}
\caption{Annual market returns from 1927 to 2023.}
\label{fig:sp500_historical}
\end{figure}

We use annual market returns from 1927 to 2023 obtained from Kenneth French's
website~\cite{french_data_lib}. 
These market returns are shown in
figure~\ref{fig:sp500_historical}, with a few statistics given in the first row of
table~\ref{tab:summary_stats}.
\paragraph{Model.}
We model the market returns as a three-component Gaussian mixture model
(GMM),
fit using Scikit-learn~\cite{scikit-GMM}. This model allows us to sample
realistic market returns.
The fitted GMM component means are 28\%, -11\%, and 11\%, and the 
corresponding standard deviations are 11\%, 16\%, and 12\%. The component
weights are 0.38, 0.25, and 0.38, respectively.

\begin{figure}
\centering
\includegraphics[width=0.6\textwidth]{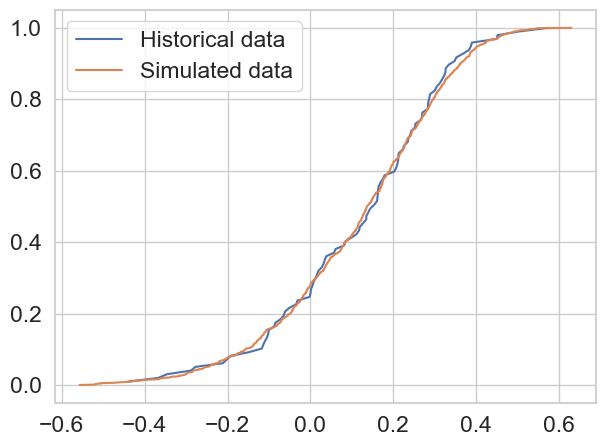}
\caption{CDF of historical and simulated market returns.}
\label{fig:cdf}
\end{figure}

\paragraph{Simulation.}
Using the GMM we can simulate realistic market returns.
A few statistics for 1,000 simulated returns are given in the second row of
table~\ref{tab:summary_stats}.
\begin{table}
\centering
\begin{tabular}{lllccccc}
\toprule
\multicolumn{3}{c}{} & \multicolumn{5}{c}{\textbf{Percentiles}} \\
\cmidrule(lr){4-8}
\textbf{Data} & \textbf{Mean} & \textbf{Volatility} & \textbf{10th} & \textbf{30th} & \textbf{50th} & \textbf{70th} & \textbf{90th} \\
\midrule
Historical & 12.0\% & 20.0\% & -11.5\% & 1.5\% & 15.4\% & 24.2\% & 33.9\% \\
Simulated & 11.7\% & 20.4\% & -15.7\% & 1.4\% & 13.6\% & 24.2\% & 36.3\% \\
\bottomrule
\end{tabular}
\caption{Market return statistics.}
\label{tab:summary_stats}
\end{table}
As seen, these statistics are reasonably close to the historical statistics.
We show the cumulative distribution
function (CDF) of the historical and simulated returns in figure~\ref{fig:cdf}.
Here too we see a good match.

\paragraph{Forecast.} It is extremely challenging to forecast market returns,
so we rely on the historical average return of 12.0\% as the (constant)
forecast for the market return.
As mentioned above, we can use a more conservative estimate, such as the 
40th percentile.

\subsection{Treasury rates and inflation}  \label{sec:treasury_inflation}   
\begin{figure}
\centering
\includegraphics[width=0.6\textwidth]{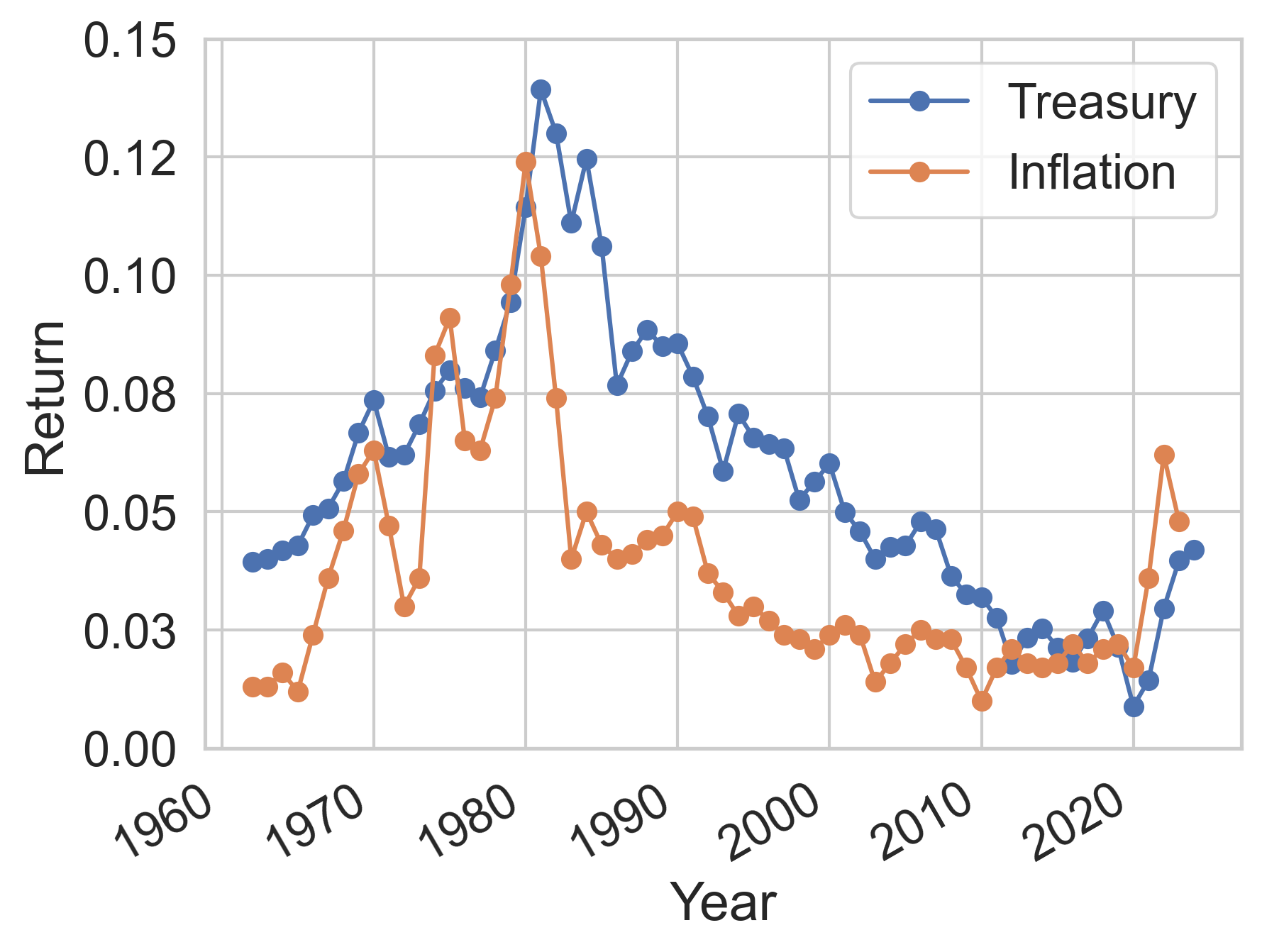}
\caption{Historical 10-year Treasury and inflation rates.}
\label{fig:treasury_inflation}
\end{figure}

We obtain 10-year Treasury rates from
1962 to 2023 from Yahoo Finance~\cite{yahoo_finance},
and annual inflation rates over the same period from the Bureau of Labor
Statistics~\cite{bls2024}. 
Figure~\ref{fig:treasury_inflation} shows the historical 10-year Treasury and
inflation rates from 1962 to 2023. 
The top rows of tables~\ref{tab:treasury_stats} and~\ref{tab:inflation_stats} show a
few statistics for the Treasury and inflation rates.
\begin{table}
\centering
\begin{tabular}{lllccc}
\toprule
\multicolumn{3}{c}{} & \multicolumn{3}{c}{\textbf{Percentiles}} \\
\cmidrule(lr){4-6}
\textbf{Data} & \textbf{Mean} & \textbf{Volatility} & \textbf{25th} & \textbf{50th} & \textbf{75th} \\
\midrule
Historical & 5.9\% & 3.0\% & 4.0\% & 5.6\% & 7.6\% \\
Simulated  & 5.3\% & 2.8\% & 3.4\% & 5.2\% & 7.2\% \\
\bottomrule
\end{tabular}
\caption{Treasury return statistics.}
\label{tab:treasury_stats}
\end{table}

\begin{table}
\centering
\begin{tabular}{lllccc}
\toprule
\multicolumn{3}{c}{} & \multicolumn{3}{c}{\textbf{Percentiles}} \\
\cmidrule(lr){4-6}
\textbf{Data} & \textbf{Mean} & \textbf{Volatility} & \textbf{25th} & \textbf{50th} & \textbf{75th} \\
\midrule
Historical & 3.8\% & 2.5\% & 2.1\% & 2.9\% & 4.8\% \\
Simulated  & 3.5\% & 2.2\% & 2.0\% & 2.6\% & 4.8\% \\
\bottomrule
\end{tabular}
\caption{Inflation rate statistics.}
\label{tab:inflation_stats}
\end{table}

\begin{figure}
\centering
\begin{subfigure}{0.48\textwidth}
\includegraphics[width=0.9\textwidth]{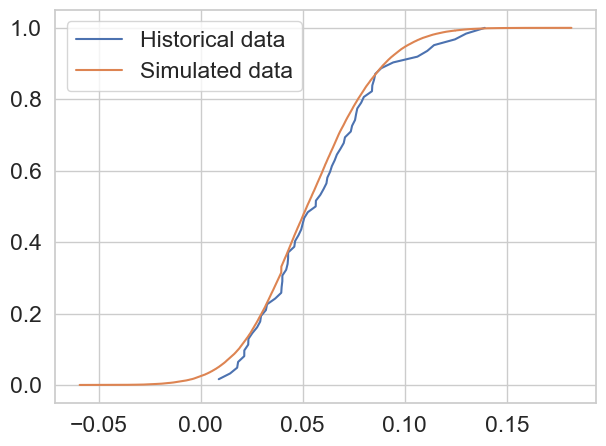}
\caption{Treasury rates.}
\end{subfigure}
\begin{subfigure}{0.48\textwidth}
\includegraphics[width=0.9\textwidth]{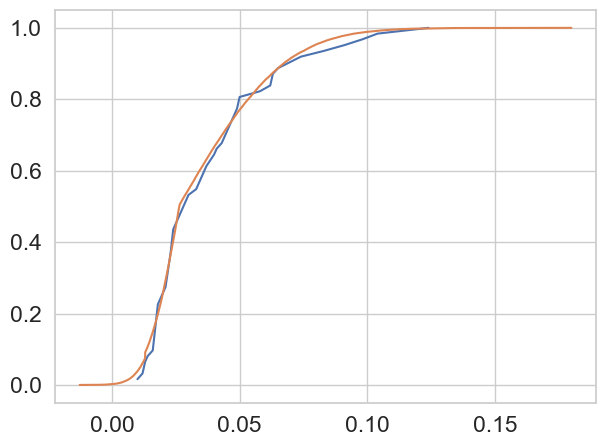}
\caption{Inflation rates.}
\end{subfigure}
\caption{CDF of historical and simulated Treasury and inflation rates.}
\label{fig:cdf_treasury_inflation}
\end{figure}

\paragraph{Model.}
The distribution of inflation rates is substantially nonsymmetric around its
median value, as can be seen in the righthand plot 
in figure \ref{fig:cdf_treasury_inflation}.
To mitigate this we first
transform the inflation rates using a simple piecewise linear 
(PWL) function $\varphi$ with
different slopes below and above the median,
\BEQ \label{eq:inflation_transform}
\varphi(z) = \left\{
\begin{array}{ll}
(z-k)s_- + k & z \leq k, \\
(z-k)s_+ + k & z > k,
\end{array}
\right.
\EEQ
where $k$ (the `kink point') is the empirical median $k=2.9\%$,
and $s_-$ and $s_+$ are the slopes of the PWL
function below and above the median, which we took as
$s_-=2.5$ and $s_+=0.75$. These values were chosen to make the transformed
inflation rates approximately, or at least closer to, Gaussian.

We model the Treasury rates and
transformed inflation rates, denoted as the vector $x_t \in \reals^2$,
as a first-order vector autoregressive (VAR) process that includes a mean.
This model has the form
\[
x_{t+1} = \mu + A (x_t-\mu) + \epsilon_t,
\]
where $\mu\in \reals^2$ is the mean, $A\in \reals^{2 \times 2}$
is the VAR (matrix) coefficient, 
and $\epsilon_t \sim \mathcal N(0,\Sigma^\epsilon)$,
where $\Sigma^\epsilon$ is the covariance matrix of the residuals.
We take $\mu$ to be the empirical mean of $x_t$,
and $A$ is chosen to minimize the sum of squared errors on the historical data.
We take $\Sigma^\epsilon$ to be the empirical covariance of the
residuals of the VAR model.  The fitted parameter values are
\[
\mu = \left[\begin{array}{c} 0.058 \\ 0.029  \end{array}\right], 
\qquad A = \left[\begin{array}{rr} 0.80 & 0.24  \\ -0.04 & 0.88 \end{array}\right],
\qquad 
\Sigma^{\epsilon} = 10^{-4} \left[\begin{array}{cc}
0.72 & 0.48 \\
0.48 & 1.47
\end{array}
\right],
\]
respectively. 

\paragraph{Simulation.}
Given initial values for the Treasury and inflation rates, we can simulate the
rates over time. We first transform the initial inflation rate using the PWL
function $\varphi$ given in \eqref{eq:inflation_transform},
then simulate the Treasury and transformed inflation rates using the VAR
model. To obtain the simulated inflation rates, we apply the inverse PWL function 
$\varphi^{-1}$ to the second component of $x_t$.

\begin{figure}
\centering
\begin{subfigure}{0.48\textwidth}
\includegraphics[width=0.9\textwidth]{figures/treasury_inflation.png}
\caption{Historical rates.}
\end{subfigure} \\
\begin{subfigure}{0.48\textwidth}
\includegraphics[width=0.9\textwidth]{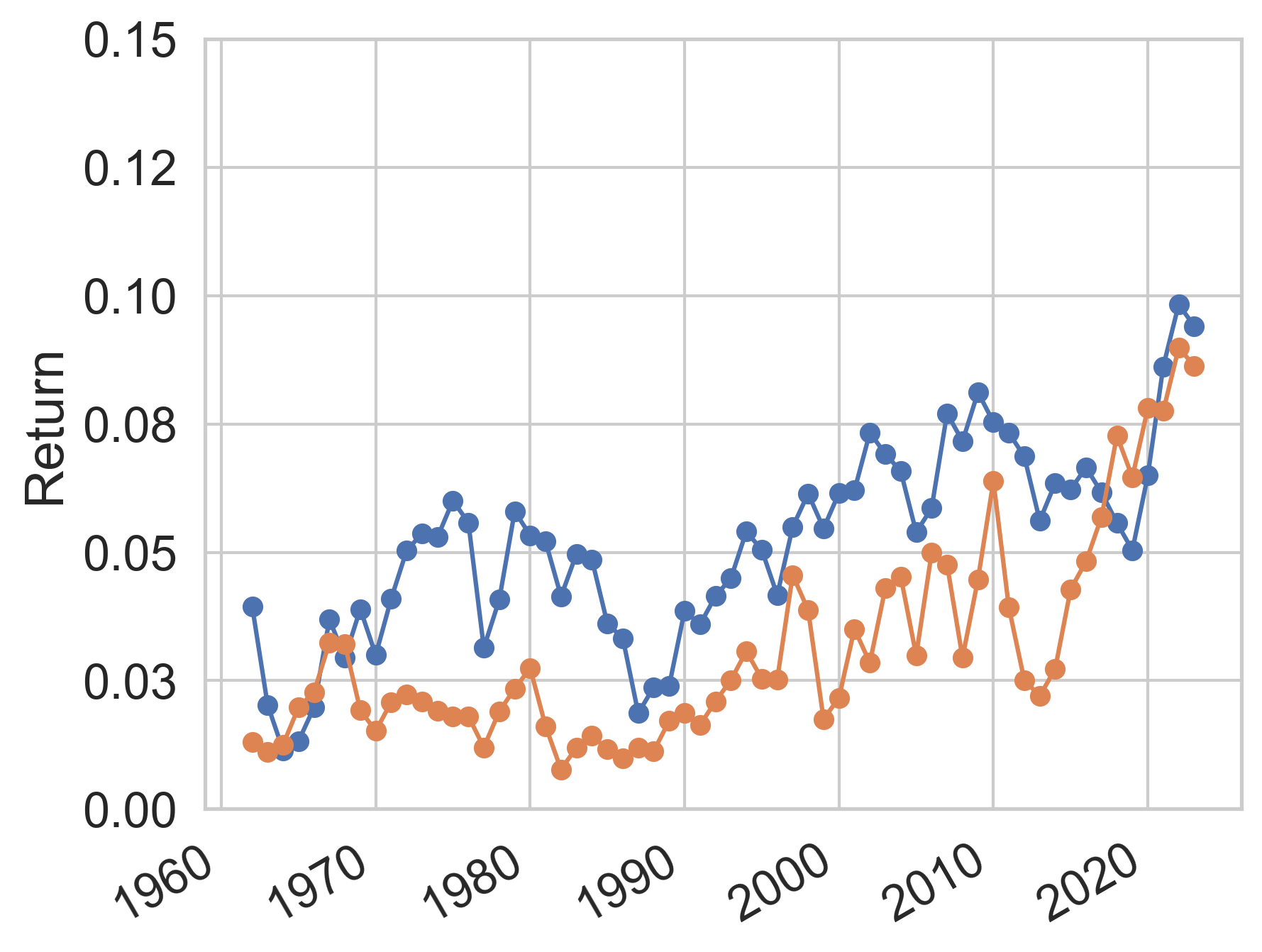}
\caption{First simulated trajectory.}
\end{subfigure} \\
\begin{subfigure}{0.48\textwidth}
\includegraphics[width=0.9\textwidth]{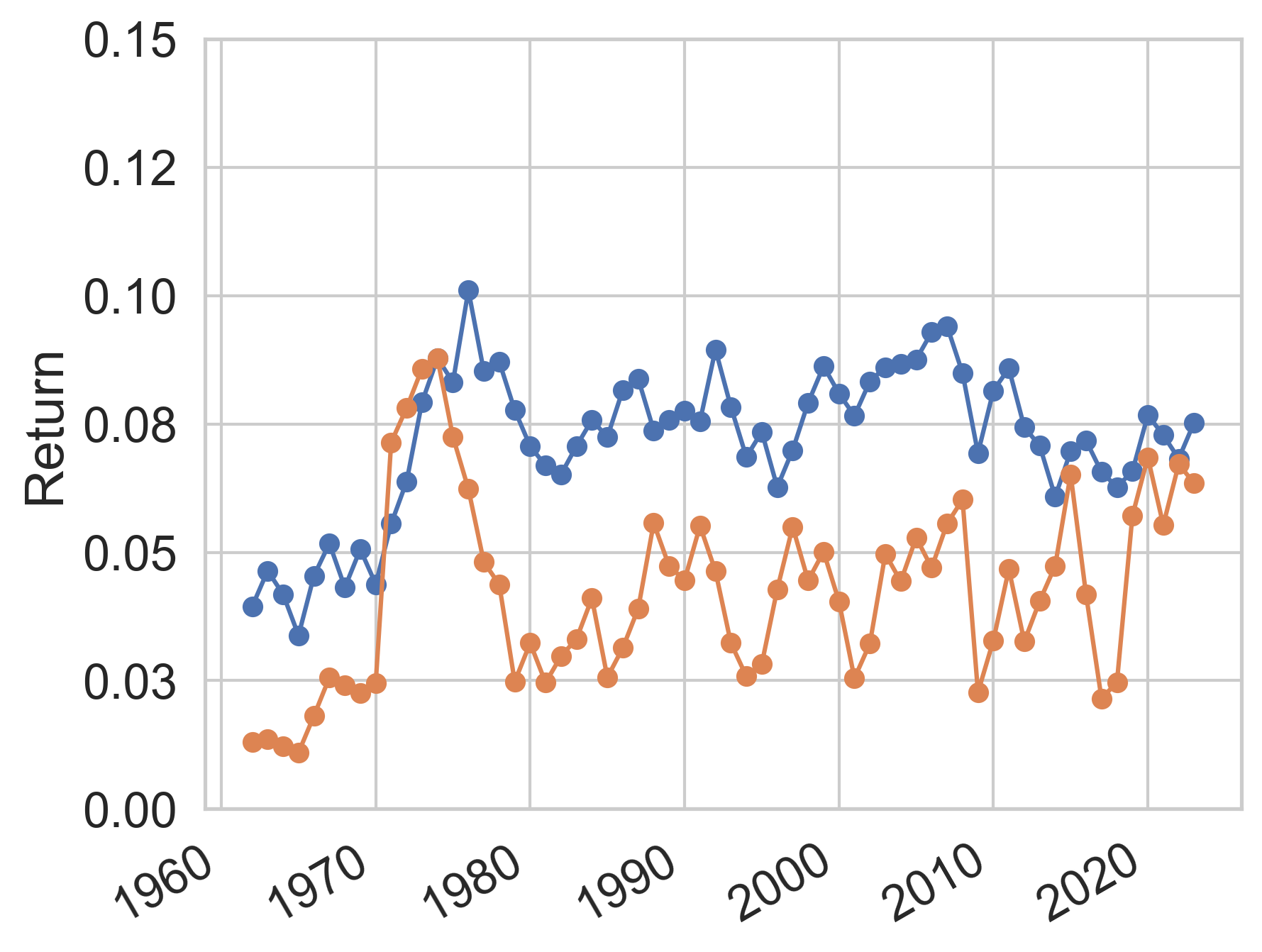}
\caption{Second simulated trajectory.}
\end{subfigure}
\caption{Two simulated trajectories of Treasury and inflation rates, initialized
from the 1962 values.}
\label{fig:treasury_inflation_sim}
\end{figure}

Two examples of simulated trajectories, initialized at the 1962
values, are shown in figure~\ref{fig:treasury_inflation_sim}. A few statistics
for 1,000 simulated trajectories, initialized at the 1962 values and simulated
until 2023, are given in tables~\ref{tab:treasury_stats}
and~\ref{tab:inflation_stats}. 
Figure~\ref{fig:cdf_treasury_inflation} shows the CDFs of the historical and
simulated rates. We can see good fits for both the Treasury and inflation rates.
As an example of a joint statistic, the correlation between the Treasury and
inflation rates is 72\% for the historical data and 70\% for the simulated data.

When simulating trajectories not initialized from specific known values,
we can sample initial Treasury and
inflation rates from their steady state distribution under 
the VAR model, which is 
$\mathcal{N}(\mu, \Sigma^\text{ss})$,
where $\Sigma^\text{ss}$ solves the Lyapunov equation
$\Sigma^\text{ss} = A \Sigma^\text{ss} A^T + \Sigma^{\epsilon}$.

\begin{figure}
\centering
\begin{subfigure}{0.48\textwidth}
\includegraphics[width=0.9\textwidth]{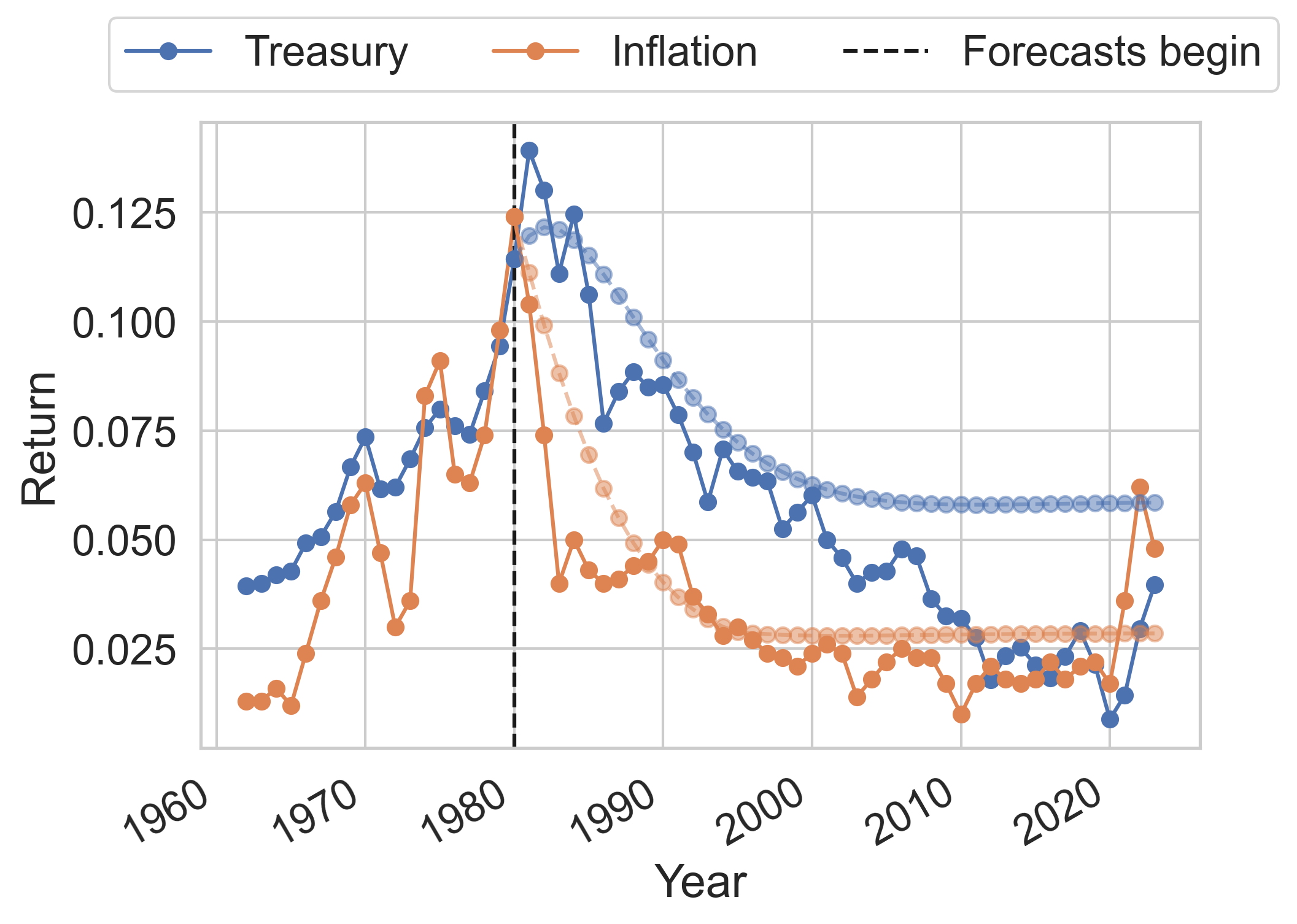}
\caption{Forecasted rates in year 1980.}
\end{subfigure}
\begin{subfigure}{0.45\textwidth}
\includegraphics[width=0.9\textwidth]{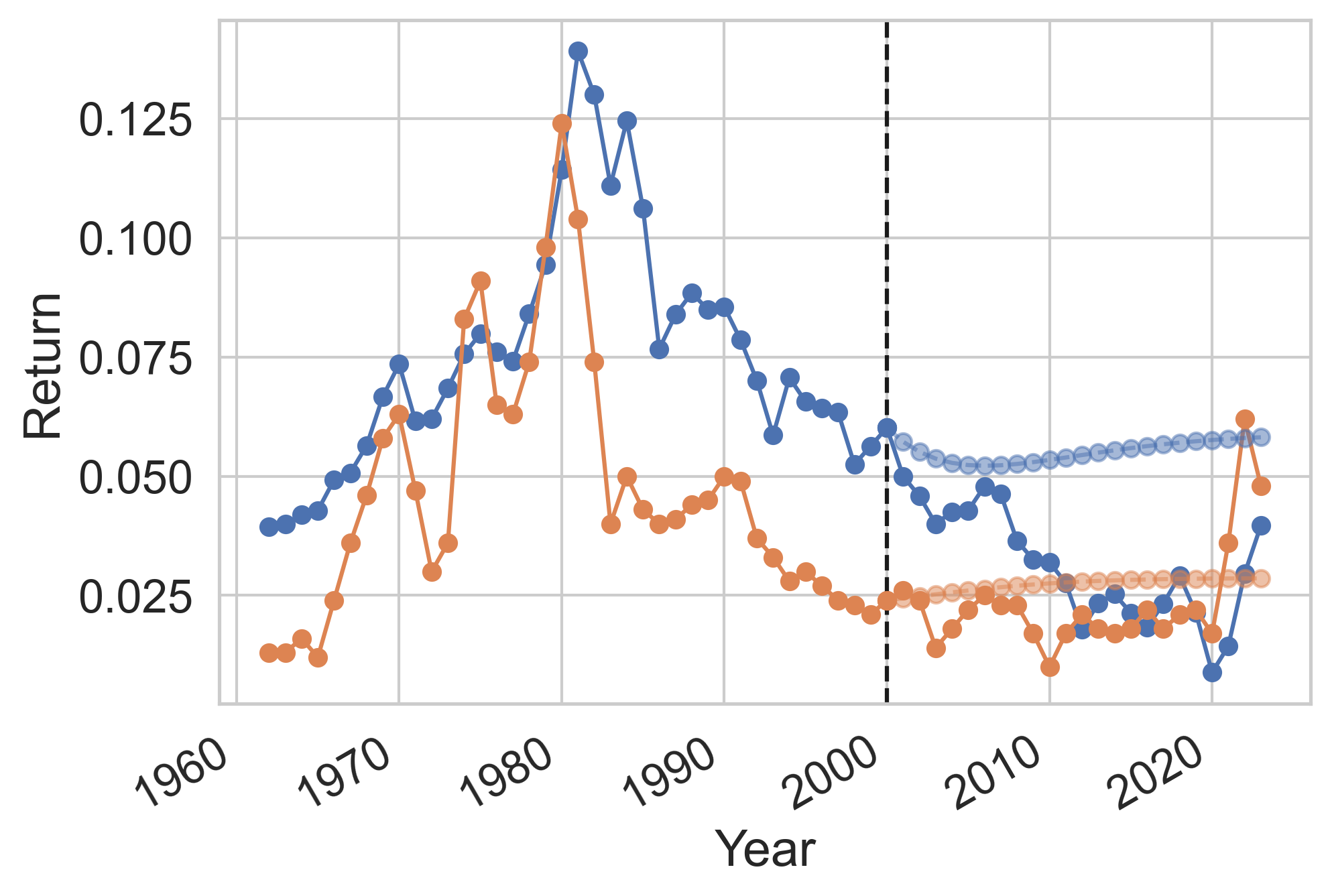}
\caption{Forecasted rates in year 2000.}
\end{subfigure}
\caption{Forecasted Treasury and inflation rates. The
shaded and solid curves represent forecasted and realized rates, respectively.}
\label{fig:forecast}
\end{figure}

\paragraph{Forecast.}
We can use the VAR model to give an elementary forecast of future Treasury and 
inflation rates. Denoting our estimate of $x_\tau$ in time period $t$, 
with $\tau>t$, as $\hat x_{\tau|t}$, we have
\[
\hat x_{\tau|t} = \mu + A^{t-\tau} (x_t-\mu), \quad \tau= t+1, t+2, \ldots.
\]
(To get the inflation rate prediction we apply $\varphi^{-1}$ 
to the second component of $\hat x_{\tau|t}$.)
Figure \ref{fig:forecast} shows the forecasted
Treasury and inflation rates, up to 2023, at two periods in time,
1980 and 2000.
They are not particularly good forecasts, but they do capture the general shape
of future values, which is all that is needed in MPC.

\subsection{Portfolio returns} \label{sec:portfolio_returns}

\begin{figure}
\centering
\includegraphics[width=0.6\textwidth]{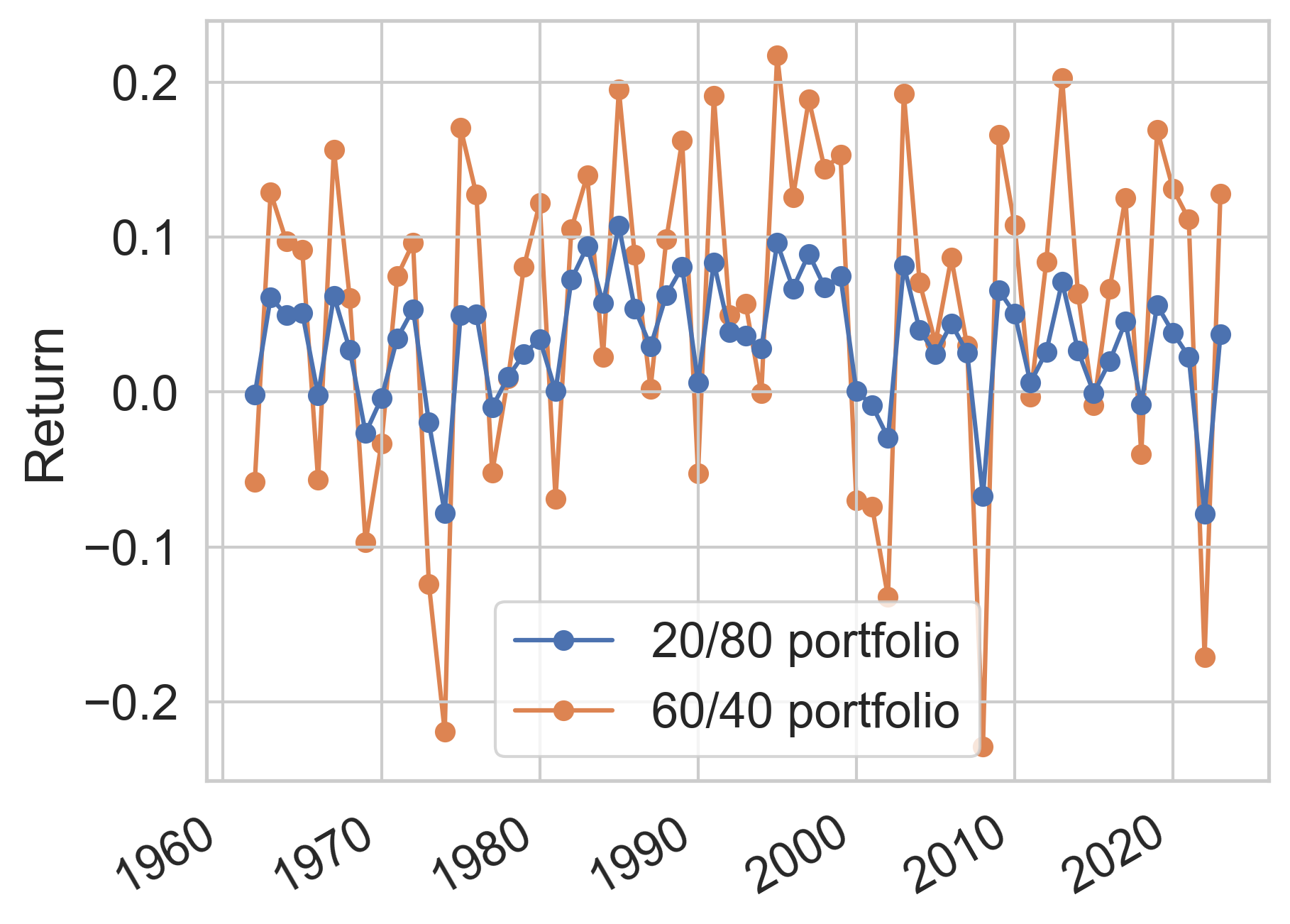}
\caption{Inflation adjusted portfolio returns.}
\label{fig:portfolio_returns}
\end{figure}

Using the statistical models for market returns, Treasury rates, and inflation
we construct models for two portfolios: 20/80 and 60/40 portfolios, comprised of
20\% and 60\% stocks, respectively, and the remainder in Treasury bonds. The
inflation adjusted returns for the portfolios is the weighted sum of the market
and Treasury returns minus inflation. The realized returns for the three
portfolios are shown in figure~\ref{fig:portfolio_returns}, with a few
statistics given in the top rows of tables~\ref{tab:2080_portfolio_stats}
and~\ref{tab:6040_portfolio_stats}, respectively.

\begin{figure}
\centering
\begin{subfigure}{0.48\textwidth}
\includegraphics[width=0.9\textwidth]{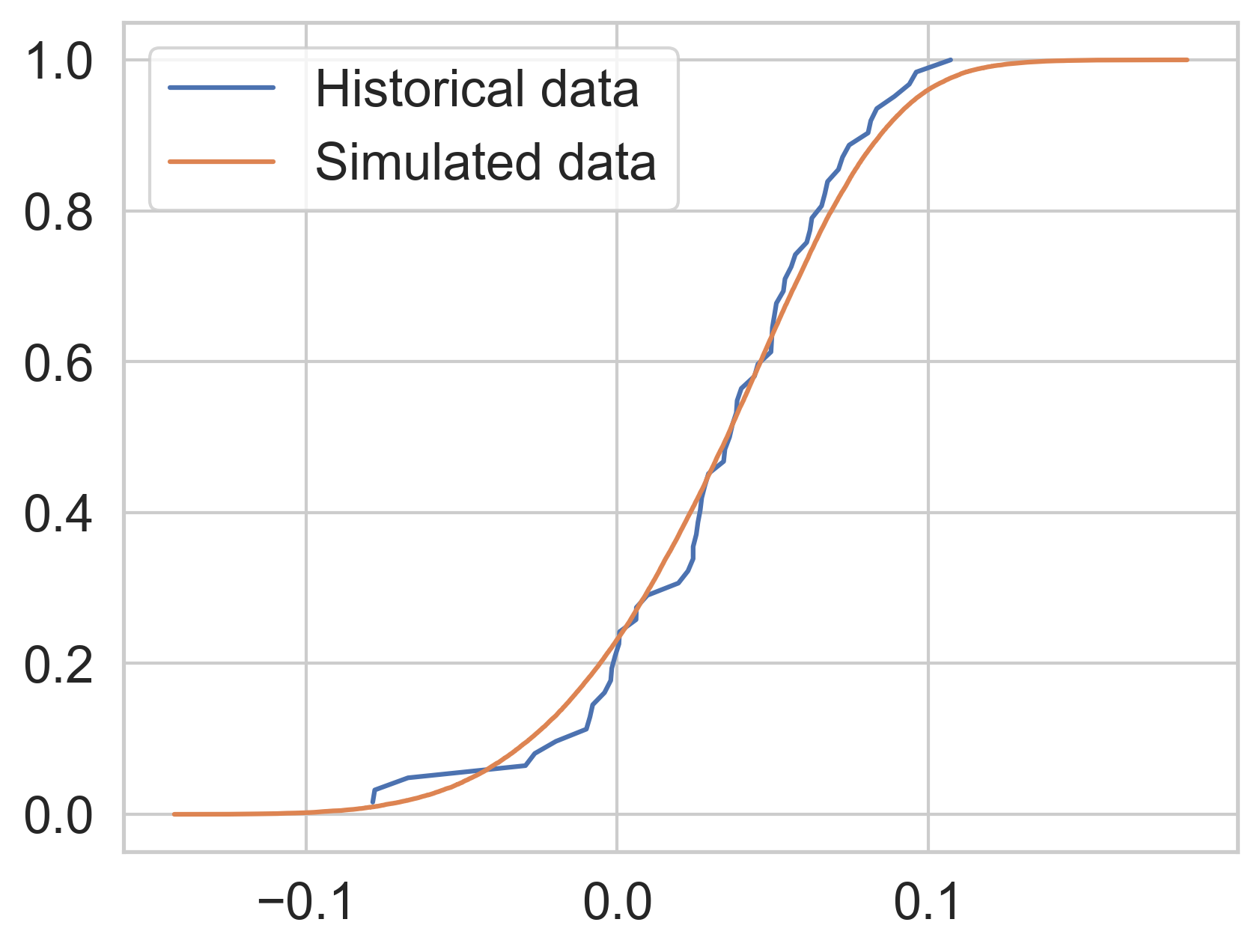}
\caption{20/80 portfolio.}
\end{subfigure} 
\begin{subfigure}{0.48\textwidth}
\includegraphics[width=0.9\textwidth]{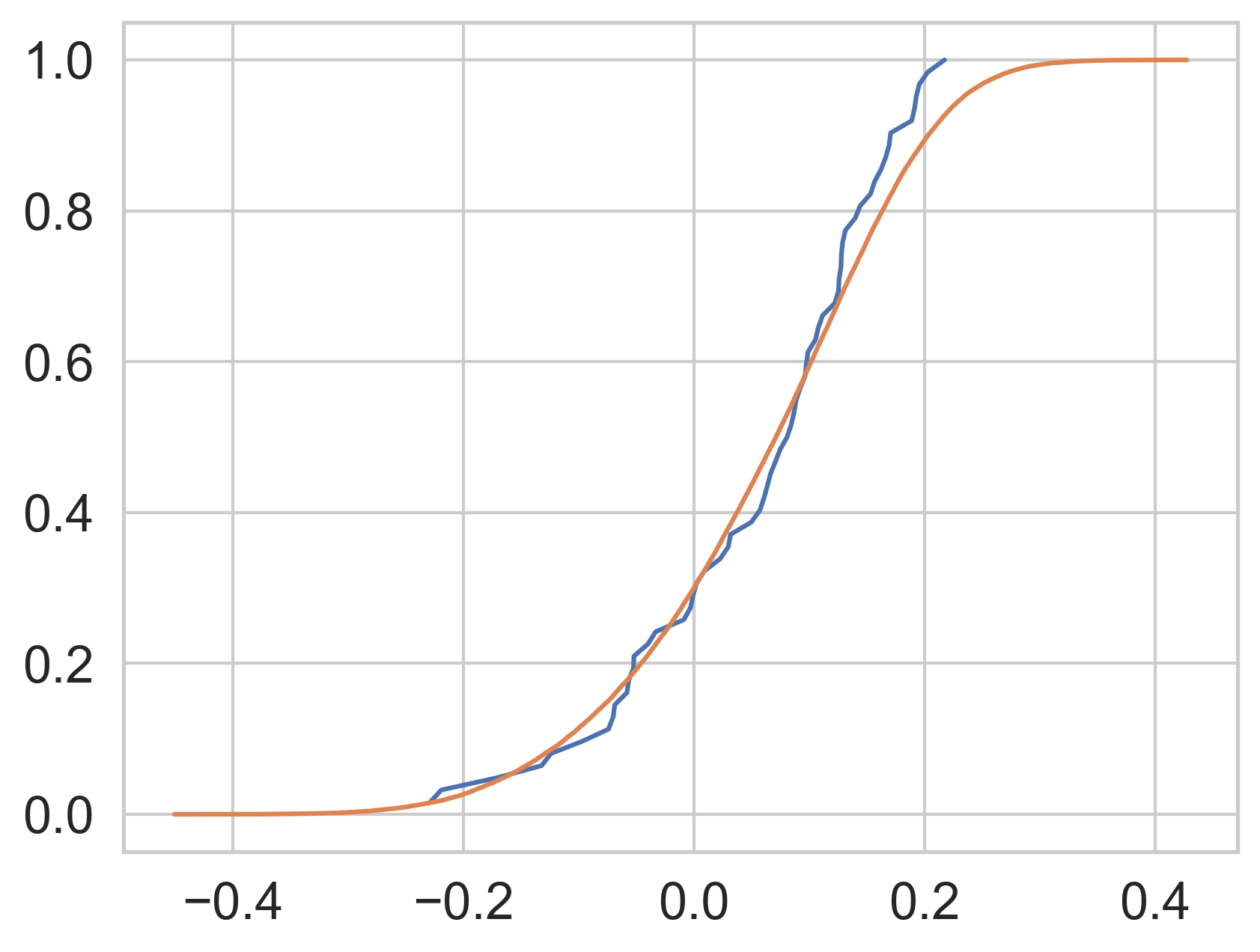}
\caption{60/40 portfolio.}
\end{subfigure}
\caption{CDF of inflation adjusted portfolio returns.}
\label{fig:cdf_portfolio_returns}
\end{figure}

\paragraph{Simulation.} We can simulate the portfolio returns by sampling from
the GMM for the market returns and the VAR model for the Treasury and inflation
rates. The CDFs of simulated and realized inflation adjusted portfolio returns
are shown in figure~\ref{fig:cdf_portfolio_returns}, and some statistics are
given in tables~\ref{tab:2080_portfolio_stats}
and~\ref{tab:6040_portfolio_stats}. (We ran 1,000 simulations initialized at the
1962 values for the Treasury and inflation rates and simulated until 2023.)

Note that the market return model is fit on data from
1927, while the Treasury and inflation rate model is fit on data from 1962.
Realized portfolio returns are available only from 1962, accounting for the
small discrepancy between the historical and simulated portfolio returns. (The
historical portfolio returns are available from 1962, while the simulated
returns leverage the market return model with data from 1927.)

\paragraph{Forecast.} 
We forecast the portfolio returns by combining the forecasts for the market
returns, Treasury rates, and inflation described above. We forecast the market
using the historical average market return of 12.0\%, and use the VAR model forecasts from
\S\ref{sec:treasury_inflation} for Treasury bonds and inflation rates. The
portfolio returns are then forecasted as a weighted sum of the market and
Treasury
returns minus inflation. Another, much simpler approach, which turns out to work
just as well, is to assume fixed returns based on historical
averages, giving 3.2\% and 5.5\% for the 20/80 and 60/40
portfolios, respectively. 
\begin{table}
\centering
\begin{tabular}{lllccc}
\toprule
\multicolumn{3}{c}{} & \multicolumn{3}{c}{\textbf{Percentiles}} \\
\cmidrule(lr){4-6}
\textbf{Portfolio} & \textbf{Mean} & \textbf{Volatility} & \textbf{25th} & \textbf{50th} & \textbf{75th} \\ 
\midrule
Historical & 3.2\% & 4.0\% & 0.6\%  & 3.7\% & 6.0\% \\ 
Simulated & 3.1\% & 4.4\% & 0.3\%  & 3.5\% & 6.3\% \\ 
\bottomrule
\end{tabular}
\caption{Inflation adjusted 20/80 portfolio return statistics.}
\label{tab:2080_portfolio_stats}
\end{table}
\begin{table}
\centering
\begin{tabular}{lllccc}
\toprule
\multicolumn{3}{c}{} & \multicolumn{3}{c}{\textbf{Percentiles}} \\
\cmidrule(lr){4-6}
\textbf{Portfolio} & \textbf{Mean} & \textbf{Volatility} & \textbf{25th} & \textbf{50th} & \textbf{75th} \\ 
\midrule
Historical & 5.5\% & 10.7\% & -0.7\% & 8.2\%  & 12.9\% \\ 
Simulated & 5.7\% & 12.1\% & -2.1\% & 7.1\%  & 14.7\% \\ 
\bottomrule
\end{tabular}
\caption{Inflation adjusted 60/40 portfolio return statistics.}
\label{tab:6040_portfolio_stats}
\end{table}



\section{Numerical experiments} \label{sec:experiments}

\subsection{Experimental setup} \label{sec:setup}
We consider two different retirees: an upper-middle-class female and a
lower-middle-class male, both retiring at age 65. 
We select their income levels, account
balances, and target consumption levels as typical,
drawing from~\cite{yahoo_net_worth_retirees, investopedia_net_worth,
census_p60_279, yahoo_income_fall}. The data and code to reproduce the
experiments are available at \url{https://github.com/cvxgrp/retirement}.

\paragraph{Investments.}
We assume both retirees invest in stock and bond portfolios
in their brokerage, IRA, and Roth accounts. The brokerage account has a mix
of 20\% stocks and 80\% Treasury bonds, and the IRA and Roth accounts have 60\%
stocks and 40\% Treasury bonds. Roughly speaking, the brokerage account is the
most conservative, reducing liquidity risk in any given year, while the IRA and
Roth are more aggressive, aiming for higher long-term returns. We described how
to forecast the returns for these investments in \S\ref{sec:portfolio_returns}.
Our simulations did not show any significant differences in the results when
using the more complicated VAR model forecasts versus the historical averages.
Hence, we adopt the simpler historical average forecasts for the portfolio
returns, giving
\[
\rho^B = 1.032, \quad \rho^I = 1.055, \quad \rho^R = 1.055.
\]

\paragraph{Planning horizon.}
To incorporate a margin of conservatism, we extend the planning horizon to 150\%
of the retiree's remaining life expectancy or until age 120, whichever occurs
first. For example, if a retiree is 65 years old with an expected lifespan of 85
years, we plan until age 95. Life expectancy estimates are obtained from the
Social Security Administration (SSA) actuarial
tables~\cite{ssa_actuarial_table}.

\paragraph{Tax rates.}
We use the federal tax brackets and rates for 2024 for single individuals, 
available on the IRS website~\cite{irs2024adjustments}, 
shown in table~\ref{tab:tax_rates_single_2024}.
We model the tax function as constant over the future, which corresponds
to the brackets being annually adjusted for inflation, but the rates in
each bracket remaining the same.
We choose different capital gains tax rates for the two retirees.

\begin{table}
\centering
\begin{tabular}{cr}
\toprule
\textbf{Marginal rate} & \textbf{Income range (\$)} \\ \midrule
10\% & up to 11,600 \\ 
12\% & 11,601 -- 47,150 \\ 
22\% & 47,151 -- 100,525 \\ 
24\% & 100,526 -- 191,950 \\ 
32\% & 191,951 -- 243,725 \\ 
35\% & 243,726 -- 609,350 \\ 
37\% & over 609,350 \\ 
\bottomrule
\end{tabular}
\caption{Federal income tax rates for single individuals in 2024.}
\label{tab:tax_rates_single_2024}
\end{table}

\subsection{Benchmark policy}
We compare the MPC policy to a benchmark policy, in which the retiree withdraws
the same amount (adjusted for inflation) each year, 
without any optimization, in analogy to the popular
4\% rule~\cite{bengen1994determining}. The retiree computes 3.75\% of the
initial balance, including all account balances and projected future earnings up
until age 85, as
the annual (post-tax) withdrawal, adjusted for inflation. 
We compute the pre-tax withdrawal that yields the target post-tax consumption
by adjusting the pre-tax amount numerically.
This roughly speaking corresponds to a 4\%
pre-tax withdrawal rate, which is a common rule of thumb for retirement
planning. We use a simple withdrawal scheme, where the retiree withdraws the
required mininum amount from the IRA account, and the remaining amount is split
between the three accounts in proportion to their remaining balances, so that 
the three accounts are depleted at the same rate.

\subsection{Example: Upper-middle-class female}
We consider a female retiree who begins receiving
Social Security payments of \$3,938 monthly at age 70; the Social Security
payments are estimated using~\cite{ssa_quickcalc} for a person who makes around
\$150,000 annually prior to retirement. 
We set the initial account balances to
\[
B^{\text{init}} = 200, \quad I^{\text{init}} = 400, \quad R^{\text{init}} = 200,
\]
in thousands of USD, 
yielding a total initial balance of \$800,000. The target consumption is set
to 3.75\% of the initial balances plus projected Social Security payments
up to age 85, giving
\[
c^\text{tar} = \text{\$58,400},
\]
corresponding roughly to a planned 30-year retirement. Finally, we let
$\gamma=500$ in the objective, chosen to give good performance with a small risk
of running out of money. We take the capital gains tax rate as $\xi =15\%$ as
this is typically the rate applicable to upper-middle-class individuals, which
our experiments also support.

\paragraph{Relative consumption and bequest.}
\begin{table}
\centering
\begin{tabular}{lllccccc}
\toprule
\multicolumn{3}{c}{} & \multicolumn{4}{c}{\textbf{Percentiles}} \\
\cmidrule(lr){4-8}
 & \textbf{Min} & \textbf{Max} & \textbf{1st}  &
\textbf{5th} & \textbf{50th}  & \textbf{95th} & \textbf{99th} \\ 
\midrule
Relative bequest & 0.42 & --- & 0.58 & 0.84 & 1.06  & 1.33 & 4.19  \\ 
\bottomrule
\end{tabular}
\caption{Relative bequest statistics for the upper-middle-class
retiree.}
\label{tab:relative_metrics_upper}
\end{table}
The objective of the retiree is to leave a large bequest, while maintaining a
stable consumption throughout retirement, with a very small risk of 
needing to consume less than the target.
To directly compare the two policies, we consider the relative 
consumption and bequest for each simulation, \ie,
the consumption and bequest for the MPC policy
divided by that of the benchmark policy.
When these numbers are larger than one, the MPC policy delivered more
to the retiree (or her heirs) than the benchmark.

The relative consumption is almost always one, meaning that the MPC policy
delivers the same consumption as the benchmark policy. It differs from one in
about 2\% of the simulations, with a minimum relative consumption of 0.78 and a
maximum of 1.25. 

Some percentiles of the relative bequest are shown in
table~\ref{tab:relative_metrics_upper}. We see that the MPC policy is more
likely to deliver a noticeably larger bequest than the benchmark policy, with a
median increase of 6\%. There is no maximum
relative bequest, since the benchmark bequest can be zero, as happens in 0.6\%
of the simulations when the retiree depletes all her investments. (The MPC
bequest is never below \$10,000.)

\begin{figure}
\centering
\includegraphics[width=0.6\textwidth]{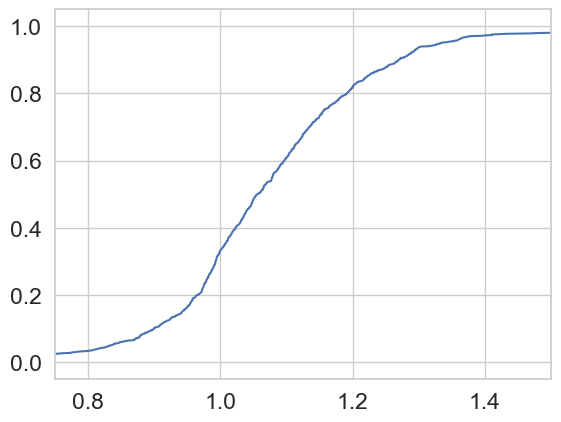}
\caption{Empirical CDF of the relative bequest for the upper-middle-class
retiree.}
\label{fig:rel_consumption_bequest}
\end{figure}

We show the empirical CDF of the relative 
bequest in figure~\ref{fig:rel_consumption_bequest}.
We see that the MPC policy delivers a larger bequest in about two thirds of 
the simulations. In the cases where a larger bequest is delivered, the median
increase is about 12\%. Overall, the MPC policy is more likely to deliver a
larger bequest than the benchmark policy, with a very low risk of delivering a
lower consumption.

\begin{figure}
\centering
\begin{subfigure}{0.48\textwidth}
\includegraphics[width=\textwidth]{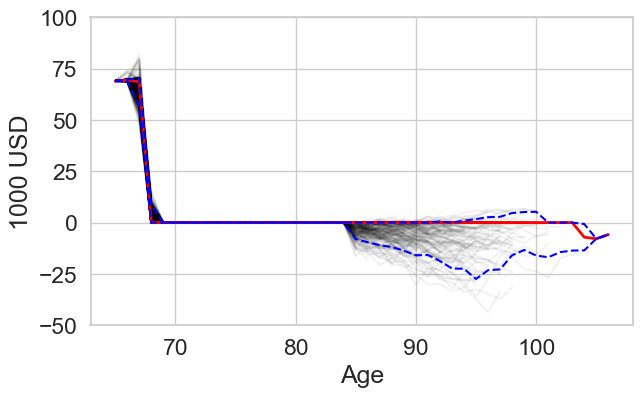}
\caption{Brokerage withdrawals, MPC policy.}
\end{subfigure}
\hspace{0.02\textwidth} 
\begin{subfigure}{0.48\textwidth}
\includegraphics[width=\textwidth]{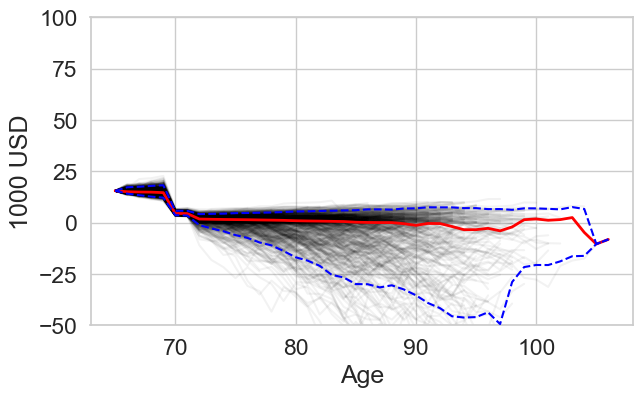}
\caption{Brokerage withdrawals, benchmark policy.}
\end{subfigure}

\vspace{0.5cm} 

\begin{subfigure}{0.48\textwidth}
\includegraphics[width=\textwidth]{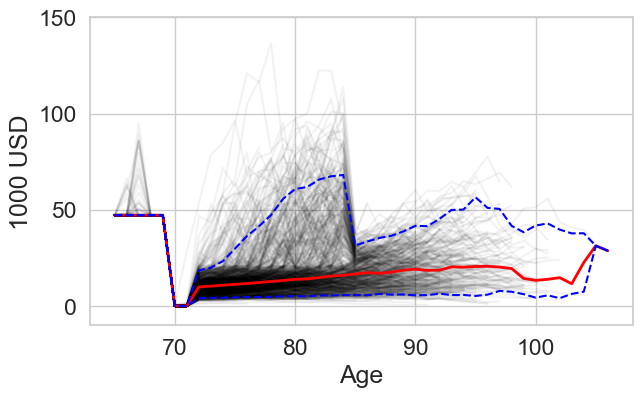}
\caption{IRA account withdrawals, MPC policy.}
\end{subfigure}
\hspace{0.02\textwidth} 
\begin{subfigure}{0.48\textwidth}
\includegraphics[width=\textwidth]{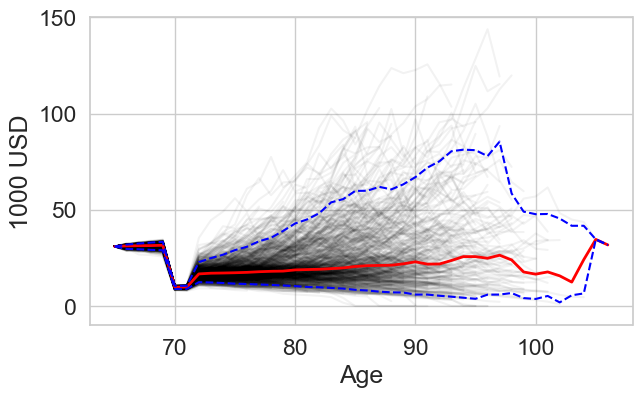}
\caption{IRA account withdrawals, benchmark policy.}
\end{subfigure}

\vspace{0.5cm} 

\begin{subfigure}{0.48\textwidth}
\includegraphics[width=\textwidth]{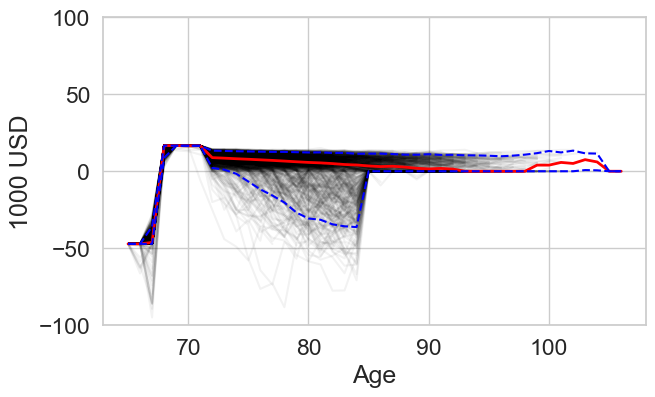}
\caption{Roth account withdrawals, MPC policy.}
\end{subfigure}
\hspace{0.02\textwidth} 
\begin{subfigure}{0.48\textwidth}
\includegraphics[width=\textwidth]{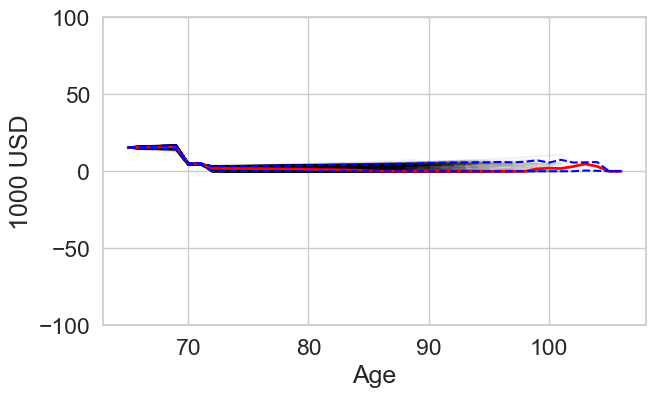}
\caption{Roth account withdrawals, benchmark policy.}
\end{subfigure}

\caption{Realized withdrawals from the brokerage, IRA, and Roth accounts for the
upper-middle-class retiree. The dark lines represent realized trajectories,
the red line represents the median over all trajectories at each age, and
the blue lines represent the 5th and 95th percentiles.}
\label{fig:withdrawals}
\end{figure}

\paragraph{Withdrawals.} The realized withdrawals for both policies from the brokerage, IRA, and
Roth accounts are shown in figure~\ref{fig:withdrawals}.
With the MPC policy, we see that the retiree withdraws from the brokerage and
IRA accounts in the first few years, until around age 70. Meanwhile, yearly deposits of up to around
\$50,000 are made to the Roth account until a little before age 70, when Roth
withdrawals start at around \$10,000--\$15,000 per year. The benchmark policy
withdrawals differ significantly from the MPC policy. No Roth conversions are
made, and rather than optimizing for tax efficiency, the retiree withdraws from
each account proportionally to the remaining balance, with any excess money
going to the brokerage account. (Of course, this is all by construction of the
benchmark policy.)

\begin{figure}
\centering
\includegraphics[width=0.6\textwidth]{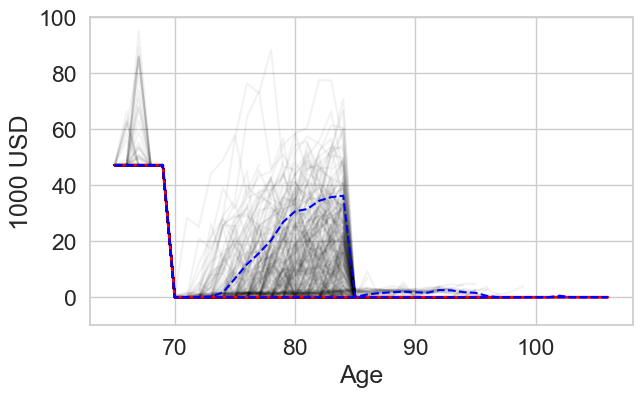}
\caption{Realized Roth conversions with the MPC policy for the
upper-middle-class retiree.
The dark lines represent realized trajectories,
the red line represents the median over all trajectories at each age, and
the blue lines represent the 5th and 95th percentiles.}
\label{fig:conversions}
\end{figure}

The MPC policy Roth conversions are shown in
figure~\ref{fig:conversions}.
Conversions start at almost \$50,000 per year and gradually decline to zero
around age 70. This pattern is intuitive: Before reaching age 70, the retiree
remains in a lower tax bracket since consumption can be funded by the brokerage
account, making conversions advantageous. However, once the brokerage account is
depleted, the retiree enters a higher tax bracket due to increased reliance on
IRA withdrawals, reducing the benefit of further conversions. After age 70, the
potential tax savings from converting no longer seem to justify the costs, given
the expected remaining lifetime, except in a few cases.

\begin{figure}
\centering
\begin{subfigure}{0.48\textwidth}
\includegraphics[width=\textwidth]{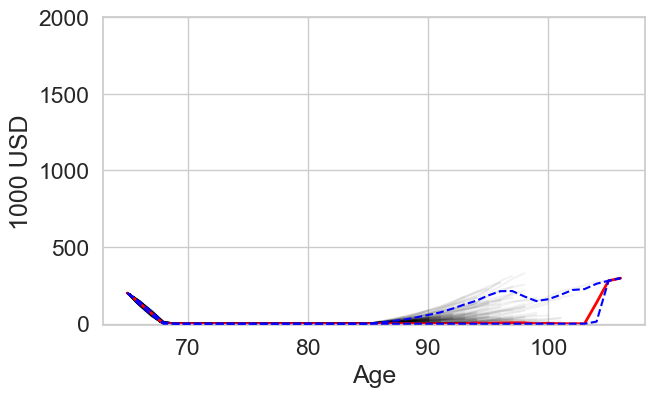}
\caption{Brokerage account balance, MPC policy.}
\end{subfigure}
\hspace{0.02\textwidth} 
\begin{subfigure}{0.48\textwidth}
\includegraphics[width=\textwidth]{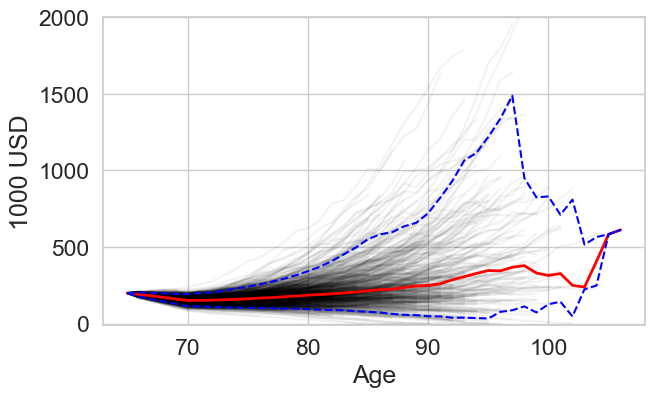}
\caption{Brokerage account balance, benchmark policy.}
\end{subfigure}

\vspace{0.5cm} 

\begin{subfigure}{0.48\textwidth}
\includegraphics[width=\textwidth]{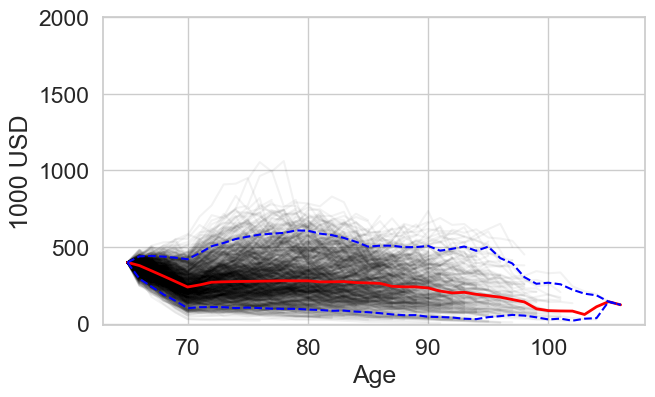}
\caption{IRA account balance, MPC policy.}
\end{subfigure}
\hspace{0.02\textwidth} 
\begin{subfigure}{0.48\textwidth}
\includegraphics[width=\textwidth]{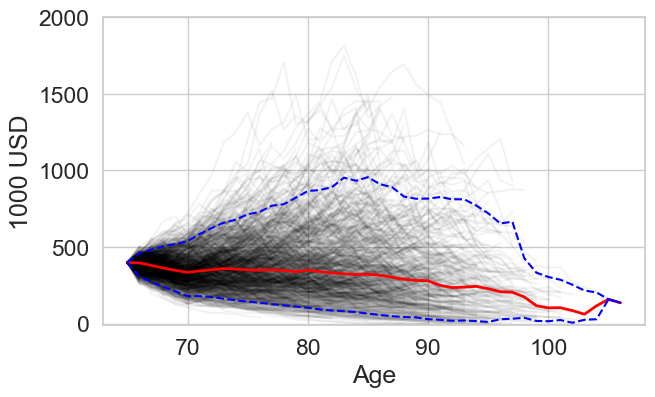}
\caption{IRA account balance, benchmark policy.}
\end{subfigure}

\vspace{0.5cm} 

\begin{subfigure}{0.48\textwidth}
\includegraphics[width=\textwidth]{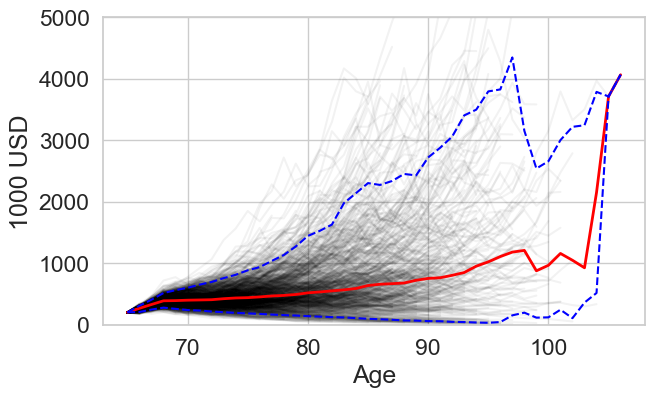}
\caption{Roth account balance, MPC policy.}
\end{subfigure}
\hspace{0.02\textwidth} 
\begin{subfigure}{0.48\textwidth}
\includegraphics[width=\textwidth]{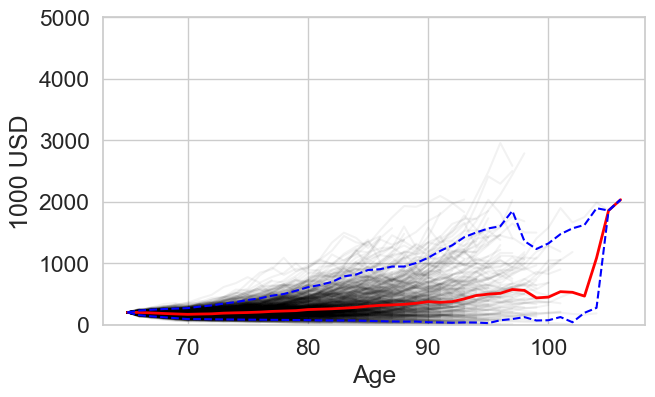}
\caption{Roth account balance, benchmark policy.}
\end{subfigure}
\caption{Realized balances of the brokerage, IRA, and Roth accounts for the
upper-middle-class retiree. The dark lines represent realized trajectories,
the red line represents the median over all trajectories at each age, and
the blue lines represent the 5th and 95th percentiles.}
\label{fig:balances}
\end{figure}

Finally, the balances of the three accounts are shown in
figure~\ref{fig:balances}.
For the MPC policy, the brokerage account is depleted around age 70, while the
IRA is significantly reduced by the same age due to substantial conversions to
the Roth account. This strategy enables the retiree to leave a large bequest,
particularly in the case of strong market returns.

\begin{figure}
\centering
\includegraphics[width=0.6\textwidth]{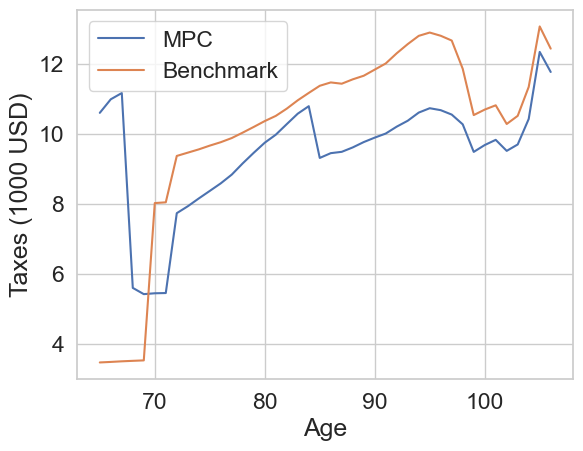}
\caption{Average realized taxes with the MPC and benchmark policies for the
upper-middle-class retiree.}
\label{fig:taxes}
\end{figure}

\paragraph{Taxes.}
Most of the improvements in the MPC policy over the benchmark policy can be
explained by the tax-efficient withdrawals and conversions. The average realized
taxes (over all 1,000 simulations) each year for the two policies are shown in
figure~\ref{fig:taxes}.
We can see that the MPC policy pays higher taxes the first few years; this is
due to the Roth conversions. After the age of 70, the MPC policy pays lower
taxes while the Roth account is allowed to grow tax-free, which is what allows
for a larger bequest.

\subsection{Example: Lower-middle-class male}
We now consider a lower-middle-class male retiree, age 65, who begins receiving Social Security
payments of \$2,013 monthly at age 70, estimated via~\cite{ssa_quickcalc} for a
person who makes \$50,000 annually prior to retiring. We let the initial account
balances be
\[
B^{\text{init}} = 50, \quad I^{\text{init}} = 100, \quad R^{\text{init}} = 0,
\]
in thousands of USD, 
yielding a total initial balance of \$150,000. The target consumption is set
to 3.75\% of the initial balance plus projected Social Security payments
up to age 85, giving
\[
c^\text{tar} = \text{\$20,100}.
\]
We let $\gamma=500$ in the objective.
We take the capital gains tax rate as $\xi = 0$, as the lower-middle-class
retiree likely is in a very low tax bracket.

\begin{figure}
\centering
\includegraphics[width=0.6\textwidth]{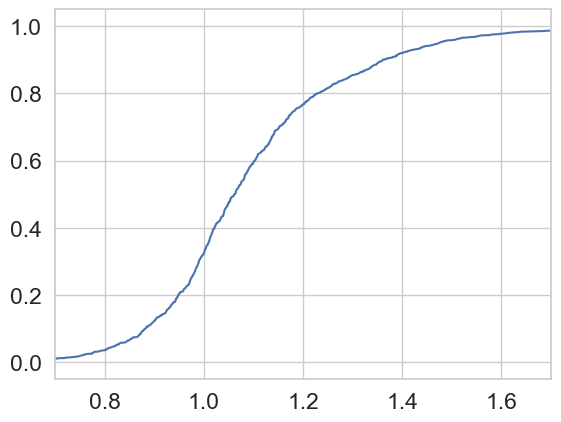}
\caption{Empirical CDF of the relative bequest for the lower-middle-class
retiree.}
\label{fig:rel_consumption_bequest_lower}
\end{figure}

\paragraph{Relative consumption and bequest.}
A few percentiles for the relative consumption and bequest of the
lower-middle-class retiree are shown in table~\ref{tab:relative_metrics_lower}.
\begin{table}
\centering
\begin{tabular}{lllccccc}
\toprule
\multicolumn{3}{c}{} & \multicolumn{4}{c}{\textbf{Percentiles}} \\
\cmidrule(lr){4-8}
 & \textbf{Min} & \textbf{Max} & \textbf{1st}  &
\textbf{5th} & \textbf{50th}  & \textbf{95th} & \textbf{99th} \\
\midrule
Relative bequest & 0.57 & 3.22 & 0.70 & 0.82 & 1.06  & 1.47 & 1.76  \\
\bottomrule
\end{tabular}
\caption{Relative bequest statistics for the lower-middle-class
retiree.}
\label{tab:relative_metrics_lower}
\end{table}
The relative consumption is almost always one, with a 0.01\% risk of going
below one. This makes sense as the lower-middle-class retiree lives off mainly the
Social Security payments alone, which he receives every year from age 70,
resulting in a
very small
risk of being forced to reduce consumption. The bequest on the other hand is higher
for the MPC policy in 68\% of the cases with a median increase of 6\%, and in
the cases where the bequest is higher, the median increase is 14\%.
The
relative bequest distribution is also heavier toward the right tail, indicating
a higher potential for a large bequest. This is further illustrated in the
empirical CDF of the relative bequest in
figure~\ref{fig:rel_consumption_bequest_lower}.

\begin{figure}
\centering
\begin{subfigure}{0.48\textwidth}
\includegraphics[width=\textwidth]{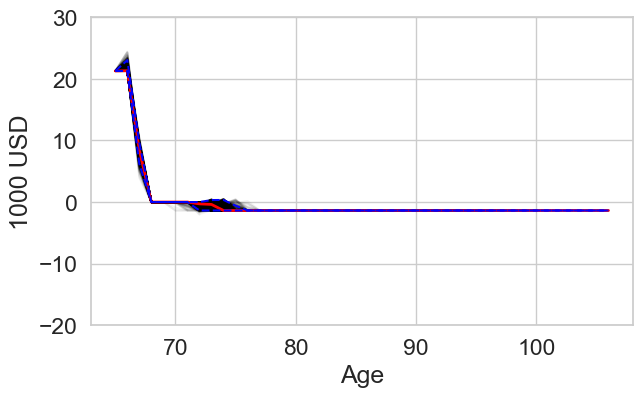}
\caption{Brokerage withdrawals, MPC policy.}
\end{subfigure}
\hspace{0.02\textwidth} 
\begin{subfigure}{0.48\textwidth}
\includegraphics[width=\textwidth]{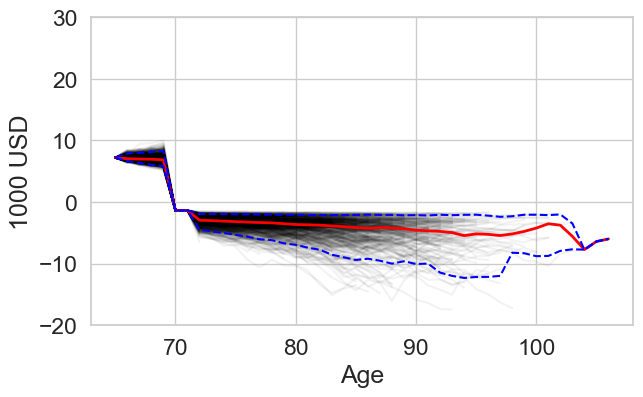}
\caption{Brokerage withdrawals, benchmark policy.}
\end{subfigure}

\vspace{0.5cm} 

\begin{subfigure}{0.48\textwidth}
\includegraphics[width=\textwidth]{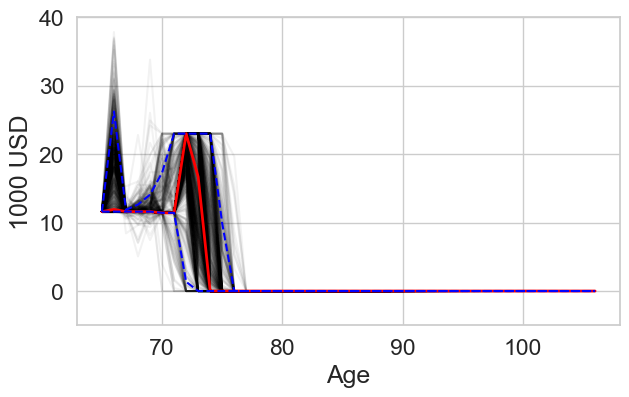}
\caption{IRA account withdrawals, MPC policy.}
\end{subfigure}
\hspace{0.02\textwidth} 
\begin{subfigure}{0.48\textwidth}
\includegraphics[width=\textwidth]{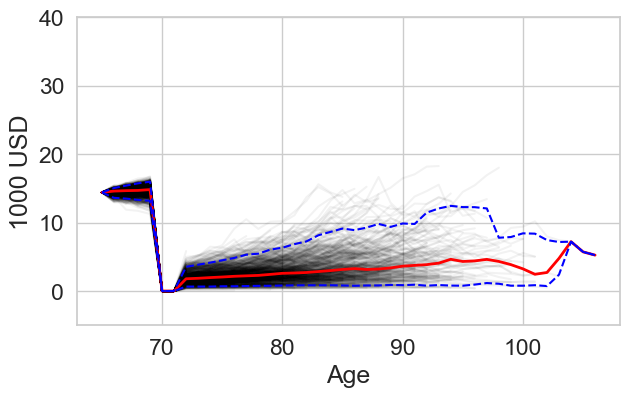}
\caption{IRA account withdrawals, benchmark policy.}
\end{subfigure}

\vspace{0.5cm} 

\begin{subfigure}{0.48\textwidth}
\includegraphics[width=\textwidth]{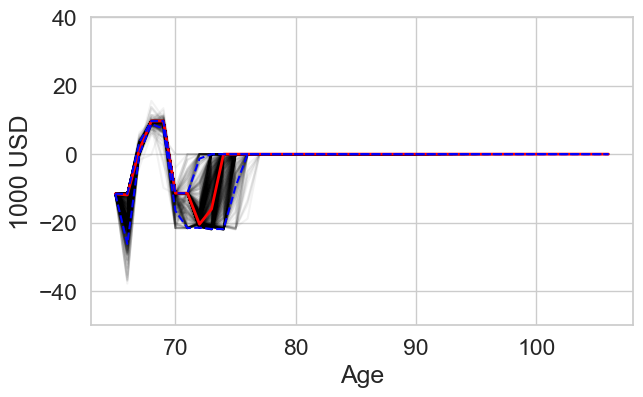}
\caption{Roth account withdrawals, MPC policy.}
\end{subfigure}
\hspace{0.02\textwidth} 
\begin{subfigure}{0.48\textwidth}
\includegraphics[width=\textwidth]{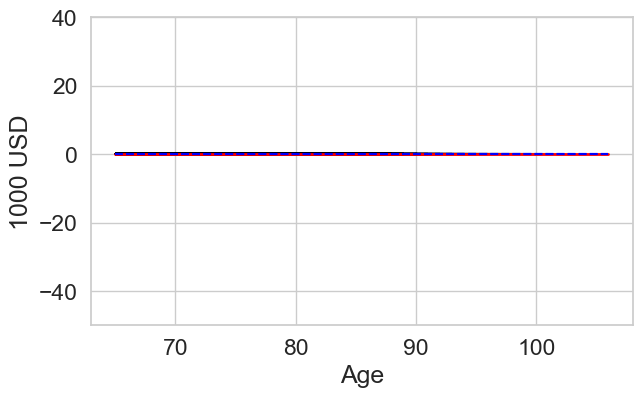}
\caption{Roth account withdrawals, benchmark policy.}
\end{subfigure}

\caption{Realized withdrawals from the brokerage, IRA, and Roth accounts for the
lower-middle-class retiree. The dark lines represent realized trajectories,
the red line represents the median over all trajectories at each age, and
the blue lines represent the 5th and 95th percentiles.}
\label{fig:withdrawals_lower}
\end{figure}

\paragraph{Withdrawals.} The realized withdrawals from the brokerage, IRA, and
Roth accounts are shown in figure~\ref{fig:withdrawals_lower}.
With the MPC policy the retiree moves the money from the IRA to the Roth account in the
first few years, until a little past age 70, via Roth conversions. From that point onwards, the retiree makes
yearly deposits to the brokerage account (with money left from the
Social Security payments after consumption), and leaves the Roth account
untouched. 
The benchmark on the other hand, withdraws from the brokerage and IRA proportionally during the first few
years, and after about age 70 makes yearly 
deposits to the brokerage account with money left from the Social
Security payments and RMDs after consumption.

\begin{figure}
\centering
\includegraphics[width=0.6\textwidth]{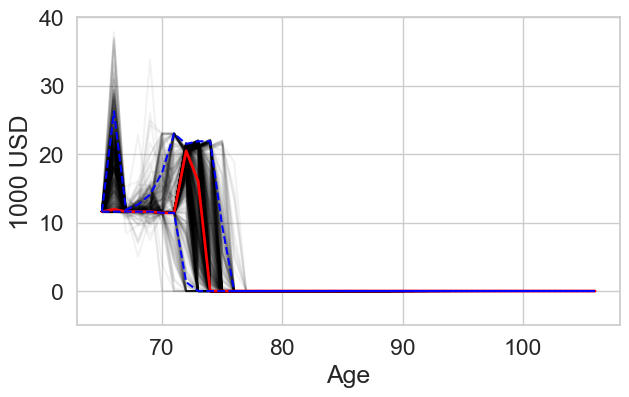}
\caption{Realized Roth conversions with the MPC policy for the
lower-middle-class retiree. The dark lines represent realized trajectories,
the red line represents the median over all trajectories at each age, and
the blue lines represent the 5th and 95th percentiles.}
\label{fig:conversions_lower}
\end{figure}

The Roth conversions are shown in
figure~\ref{fig:conversions_lower}.
They are equal to the IRA withdrawals in figure~\ref{fig:withdrawals_lower}. The
difference between the Roth withdrawals (which are negative) in
figure~\ref{fig:withdrawals_lower}
and the Roth conversions in figure~\ref{fig:conversions_lower} is the consumption.

\begin{figure}
\centering
\begin{subfigure}{0.48\textwidth}
\includegraphics[width=\textwidth]{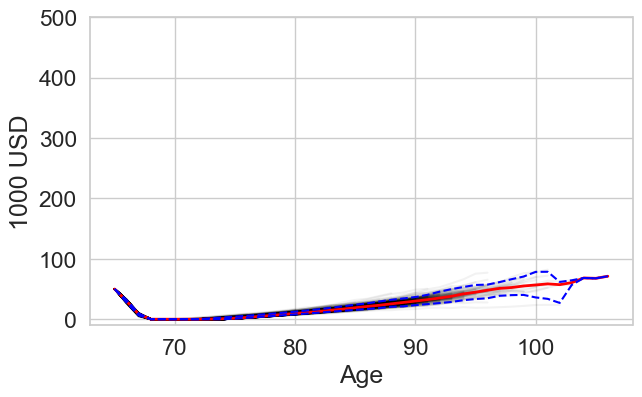}
\caption{Brokerage account balance, MPC policy.}
\end{subfigure}
\hspace{0.02\textwidth} 
\begin{subfigure}{0.48\textwidth}
\includegraphics[width=\textwidth]{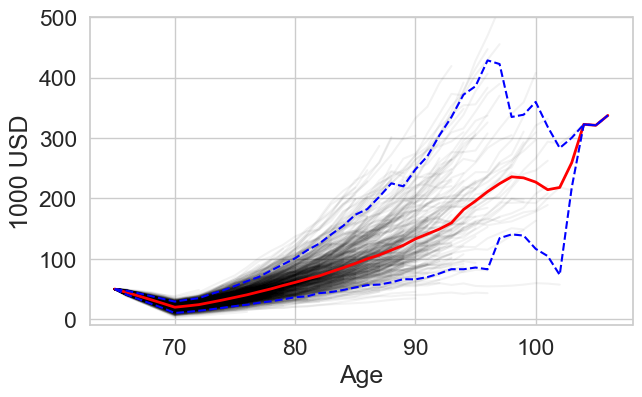}
\caption{Brokerage account balance, benchmark policy.}
\end{subfigure}

\vspace{0.5cm} 

\begin{subfigure}{0.48\textwidth}
\includegraphics[width=\textwidth]{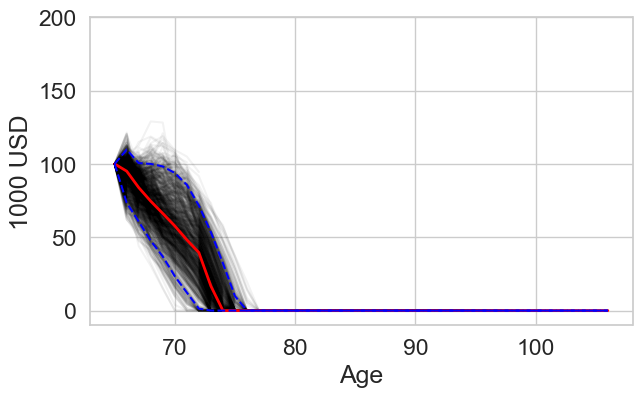}
\caption{IRA account balance, MPC policy.}
\end{subfigure}
\hspace{0.02\textwidth} 
\begin{subfigure}{0.48\textwidth}
\includegraphics[width=\textwidth]{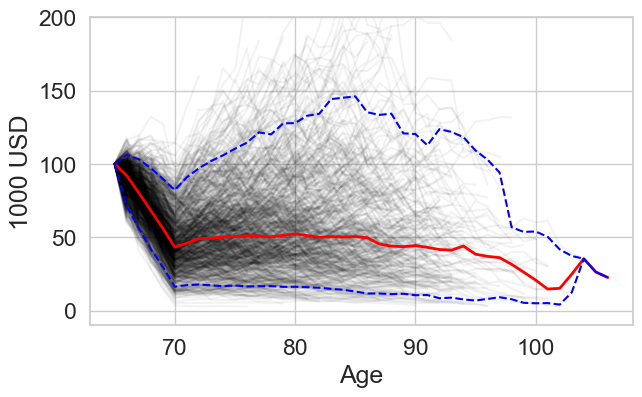}
\caption{IRA account balance, benchmark policy.}
\end{subfigure}

\vspace{0.5cm} 

\begin{subfigure}{0.48\textwidth}
\includegraphics[width=\textwidth]{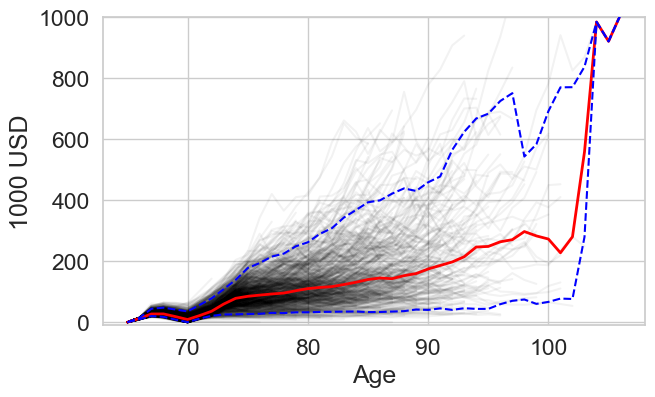}
\caption{Roth account balance, MPC policy.}
\end{subfigure}
\hspace{0.02\textwidth} 
\begin{subfigure}{0.48\textwidth}
\includegraphics[width=\textwidth]{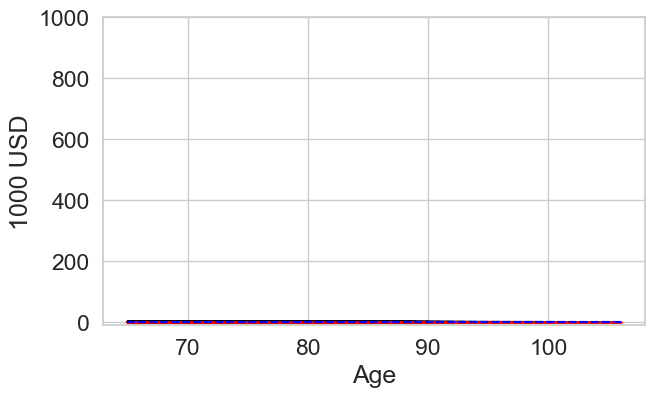}
\caption{Roth account balance, benchmark policy.}
\end{subfigure}
\caption{Realized balances of the brokerage, IRA, and Roth accounts for the
lower-middle-class retiree. The dark lines represent realized trajectories,
the red line represents the median over all trajectories at each age, and
the blue lines represent the 5th and 95th percentiles.}
\label{fig:balances_lower}
\end{figure}

Finally, the balances of the three accounts are shown in
figure~\ref{fig:balances_lower}.
With the MPC policy, the IRA account is reduced to zero a little past age 70,
increasing the Roth account balance. Once Social Security payments start, the
brokerage account starts growing, while the Roth account remains untouched,
growing with investment gains.

\begin{figure}
\centering
\includegraphics[width=0.6\textwidth]{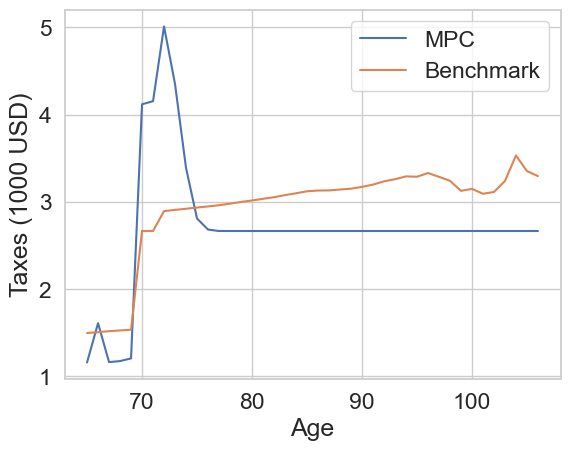}
\caption{Average realized taxes with the MPC and benchmark policies for the
lower-middle-class retiree.}
\label{fig:taxes_lower}
\end{figure}

\paragraph{Taxes.}
The average realized
taxes each year for the two policies are shown in figure~\ref{fig:taxes_lower}.
The MPC policy has higher taxes for the first few years following age 70, when a
large portion of the Roth conversions are made.
Other than that the MPC policy pays lower taxes, due to tax efficient
withdrawals.

\section{Extensions and variations}\label{sec:extensions} Here we describe some
extensions and variations of the problem, relaxing some simplifications
of problem~\eqref{prob:nonconvex}.

\subsection{Model improvements}
We first describe a few ways to improve the model in order to make it more
accurate.

\paragraph{Taxation of capital gains.}
Our model approximates capital gains taxation using a constant rate and assumes
that retirees receive no benefit from losses. In reality, the capital gains tax
rate varies with taxable income, and capital losses can be deducted against
capital gains or, in some cases, offset against ordinary income up to a
specified limit, with any excess losses carried forward to future years. As a
result, the effective tax burden on investment gains may be lower than our
approximation, depending on the retiree's financial situation. While this
simplification introduces some conservatism, incorporating more precise tax rules
into the model is relatively straightforward.

\paragraph{Five-year rules.} There are two key five-year rules: the Roth IRA
five-year rule and the Roth conversion five-year rule. The Roth IRA rule
requires that earnings can only be withdrawn tax-free if the account has been
open for at least five years, even after age 59.5. If the funds are withdrawn
before meeting this requirement, the earnings are subject to income tax and a
10\% early withdrawal penalty. The Roth conversion rule requires that converted
funds remain in the Roth account for at least five years to avoid penalties. If
withdrawn before five years, the converted amount is subject to a 10\% early
withdrawal penalty, but only for individuals under age 59.5.

We can incorporate these rules into the problem setup as follows. First, to
account for the Roth IRA rule, we add any Roth withdrawal $ r^w_t $, where $ t $
is less than five years after the Roth account was opened, to the taxable
income, and reduce $ r^w_t $ by 10\% in the cash balance equation.

Second, to account for the Roth conversion rule, we introduce two vectors. The
Roth vector $ \bar R_t \in \mathbf{R}^5 $ defines the Roth account balance in
each of the five conversion slots, tracking the balance for each year since the
money was converted. We also introduce a Roth withdrawal vector $ \bar r^w_t \in
\mathbf{R}^5 $, which contains the Roth withdrawals for each of the five
conversion slots. The total Roth account balance and  Roth withdrawal are 
\[
R_t = \mathbf{1}^T \bar R_t, \quad r^w_t = \mathbf{1}^T \bar r^w_t.
\] 
The vectors $\bar R_t$ and $\bar r^w_t$ evolve
according to linear dynamic equations. 

\paragraph{Employee sponsored accounts.}
Many retirees have additional accounts, such as a 401(k) or Roth 401(k), which
are very similar to IRAs and Roth IRAs, but have different rules, mainly higher
contribution limits and possible company matching. We can handle this by adding
two accounts in~\S\ref{sec:retirement}. They would have similar constraints to
the IRA and Roth accounts: 
\begin{itemize}
    \item The total deposit to the 401(k) and Roth 401(k) accounts is limited.
    In 2024, the limit was \$22,500 for individuals under 50, and \$30,000 for
    those aged 50 and above.
    \item Any deposit to the 401(k) and Roth 401(k) must come from earned income
    in that year.
    \item The 401(k) account is subject to RMDs starting at age 73.
    \item The Roth 401(k) account is also subject to the five-year rule for
    qualified distributions. (Roth 401(k) accounts also used to be subject to
    RMDs, unlike Roth IRAs, but this rule was recently repealed.)
\end{itemize}
We can add additional accounts and the corresponding deposit and withdrawal
constraints to the problem.

\paragraph{MAGI dependency.} In addition to the combined contribution limit for
both traditional IRA and Roth IRA accounts, there are income-based restrictions
on Roth IRA contributions for high-income earners. For 2024, an individual must
have a modified adjusted gross income (MAGI) under \$146,000 to contribute the
full amount to a Roth IRA. The contribution limit phases out for MAGI between
\$146,000 and \$161,000. Above \$161,000, no contribution is allowed. These
limits add a set of linear constraints, depending on the retiree's MAGI, which
can be estimated from the retiree's projected income and tax deductions.

There are also income-based restrictions on traditional IRA contributions. For
2024, if an individual is covered by a workplace retirement plan, they must have
a MAGI under \$66,000 to contribute the full amount to a traditional IRA. The
contribution limit phases out for MAGI between \$66,000 and \$76,000. Above
\$76,000, no tax-deductible contribution is allowed, although a non-deductible
contribution is still allowed. These IRA limits affect the taxable income in our
problem formulation for individuals who are covered by a workplace retirement
plan.

\subsection{Investments} \label{sec:investments}
In our numerical experiments we used simple stock/bond portfolios,
since these are traditional and widely used.  The few adverse
retirement outcomes we see in these simulations are all due to extreme 
left-tail events in the stock market. 
It is straightforward to intentionally construct
portfolios with less left-tail risk, by leveraging annuities, treasury
inflation protected securities (TIPS), collars, or ETFs tracking TIPS or
collar protected indexes. The methodology of this paper still applies to these
more sophisticated portfolios, with potentially better results (in terms of
worst-case outcomes) since the left-tail risk is reduced.
We present numerical results for a specific case of constructing a more
sophisticated investment portfolio in appendix \ref{sec:collar_investments}.
\paragraph{Stocks/bond split.} We could use a different mix of stocks and bonds.
For example, we could reduce the stock exposure to reduce the downside risk.
This would reduce the average bequest, but also improve the worst-case
outcomes. An alternative approach is to use age-adjusted investments. For
example, a retiree could begin with an aggressive stock/bond allocation, such as
80/20, and gradually decrease stock exposure over time. This strategy offers
greater growth potential in the early years while mitigating downside risk as
the retiree gets older. Examples of such age-adjusted investments are
Fidelity's
Freedom Funds~\cite{fidelityFreedomFunds}. 


\paragraph{Annuities.} An annuity is an insurance product that provides a fixed
stream payment (typically not inflation adjusted) over the lifetime of the retiree,
starting either at the time of purchase or at a deferred time.
Including annuities in a retirement portfolio can reduce longevity risk,
ensuring that the retiree does not outlive their savings. For example, an
immediate fixed annuity can be purchased using a portion of the retirement
savings, which guarantees a steady income regardless of market performance. 
Indexed annuities allow the retiree to participate in some of
the upside of market returns, while protecting the downside.

\paragraph{TIPS.} TIPS are U.S. government bonds designed to protect investors
from inflation risk. The principal value of TIPS increases with inflation, as
measured by the Consumer Price Index (CPI), ensuring that the purchasing power
of the investment is preserved over time. In the context of retirement planning,
incorporating TIPS into a portfolio provides a reliable, low-risk source of
income that adjusts for inflation, which is critical for retirees whose expenses
may rise over time. Additionally, ETFs that track TIPS indices can offer
liquidity and diversification benefits while maintaining the inflation
protection inherent in the underlying securities~\cite{ishares_tip}.

\paragraph{Collars.} A collar is a financial strategy that involves purchasing a
protective put option while selling a covered call option on the same underlying
asset. This strategy limits both the potential downside and upside of the
investment. In the context of retirement planning, collars can be used to
stabilize the returns on stock investments, reducing exposure to significant
market downturns. For instance, the retiree could purchase a put option limiting
the downside to -5\% annually, say, and finance this by selling a call option 
for the same price.
This limits both the downside and the upside, and like an indexed annuity 
can give returns with fewer extreme values. ETFs tracking collar-protected
indexes are also available, making it easier to gain exposure to
collar-protected
portfolios~\cite{xclr_etf}.

\section{Conclusion} \label{sec:conclusion}
We have presented a convex optimization formulation of 
the retirement funding planning problem, aimed
at maximizing a retiree's bequest while delivering a target 
inflation adjusted consumption,
accounting for taxes and rules governing the different types of 
accounts. Building on this formulation of the planning problem
we propose an MPC algorithm that dynamically
optimizes annual withdrawals and account conversions across brokerage, IRA, and
Roth accounts, with careful consideration of tax implications. 

Through extensive simulations, we demonstrated that the MPC approach
outperforms traditional withdrawal strategies, such as the 4\%
rule, in key metrics. Specifically, the MPC policy delivers comparable
consumption, with bequest outcomes that are better in around 70\% of the cases,
with a median increase around 6\%.

\subsection*{Acknowledgments}
The authors are grateful to Trevor Hastie for insightful discussions and
detailed feedback on an early draft. We thank Ron Kahn for his valuable
comments and suggestions, which greatly improved the manuscript. We also thank
Eric Luxenberg and Max Schaller for helpful discussions.

\clearpage
\bibliography{refs}

\clearpage
\appendix
\section{More sophisticated investment portfolios}
\label{sec:collar_investments}
In \S\ref{sec:investments} we described a few alternative investment
portfolios that could be used in the retirement planning problem. Here we
describe one of them in more detail, namely the collar portfolio, and show how
the results change when we use this portfolio instead of the simple
stock/bond portfolio.

\subsection{Portfolios with collar options}
A collar is an option strategy that involves buying a put option and
selling a call option on the same underlying asset, effectively limiting the
range of possible returns. The put option provides downside protection,
while the call option limits the upside potential. The downside is referred to
as
the ``floor'' ($F$) and the upside as the ``cap'' ($C$). We will assume the
investor invests in a collar on the market portfolio, \ie, the S\&P 500 index,
together with a bond (the risk-free asset).

\paragraph{Investment portfolio.} Just as before, the investor invests in
stock/bond portfolios, but she now purchases a collar on the stock portfolio.
(We do not need a collar on the bond portfolio, as the bond return is
approximately known.) We will make the collar self-financing, meaning that we
purchase a put option, limiting the downside, and sell a call option for the
same price, limiting the upside. The collar is self-financing, as the
cost of the put option is offset by the premium received from selling the call
option. For simplicity, we price the options using the Black-Scholes
formula~\cite{black_scholes_1973}.

\paragraph{Determining floor and cap.} We assume the investor wants to hedge her
stock/bond portfolio against poor market returns. She holds a weight $w$ in
stocks, and $1-w$ in bonds. She wants to ensure that the total inflation
adjusted return is at least $r_{\text{min}}$. This means that the floor $F$ must
satisfy
\BEQ \label{eq:floor}
(1-w) r_f + w F - i \geq r_{\text{min}},
\EEQ
where $r_f$ is the risk-free rate and $i$ is the inflation rate. Since we can
not have a self-financing collar with a floor greater than the risk-free rate,
we set the floor as
\BEQ \label{eq:floor2}
F = \min\left\{\frac{r_{\text{min}} - (1-w) r_f + i}{w}, r_f \right\}.
\EEQ
In words, if we can guarantee a floor greater than the risk-free rate, we do
so. Otherwise, we set the floor to the risk-free rate.

We then sell a call option with a cap $C$, where $C$ is defined by the call
option that has the same price as the put option with floor $F$. In summary we
find the floor and cap in three steps.
\BIT
\item Determine the minimum inflation adjusted return the investor is willing to
accept, $r_{\text{min}}$.
\item Find the floor $F$ using~\eqref{eq:floor2}.
\item The cap $C$ is the strike price of the call option with the same
price as the put option with floor $F$.
\EIT

\subsection{Experimental setup}
\begin{figure}
\centering
\begin{subfigure}{0.48\textwidth}
\includegraphics[width=0.9\textwidth]{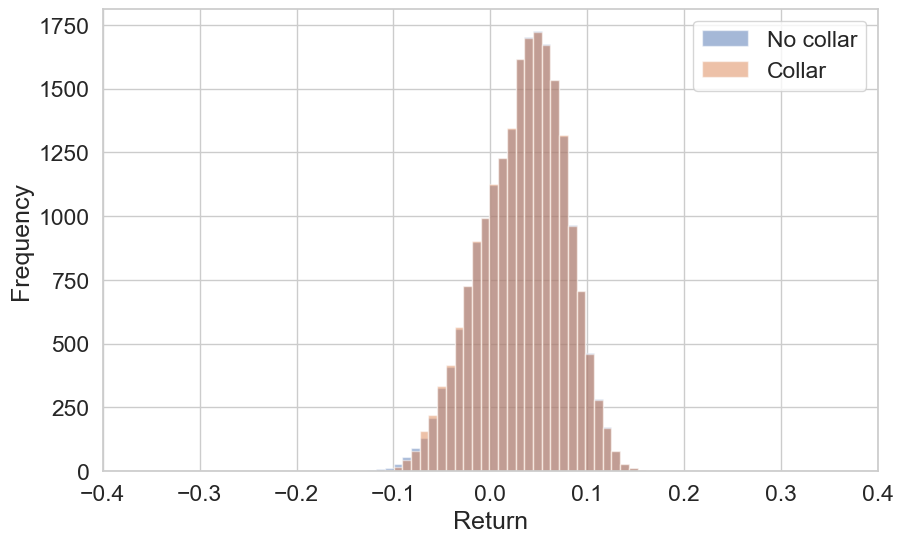}
\caption{20/80 stock/bond portfolio.}
\end{subfigure}
\begin{subfigure}{0.48\textwidth}
\includegraphics[width=0.9\textwidth]{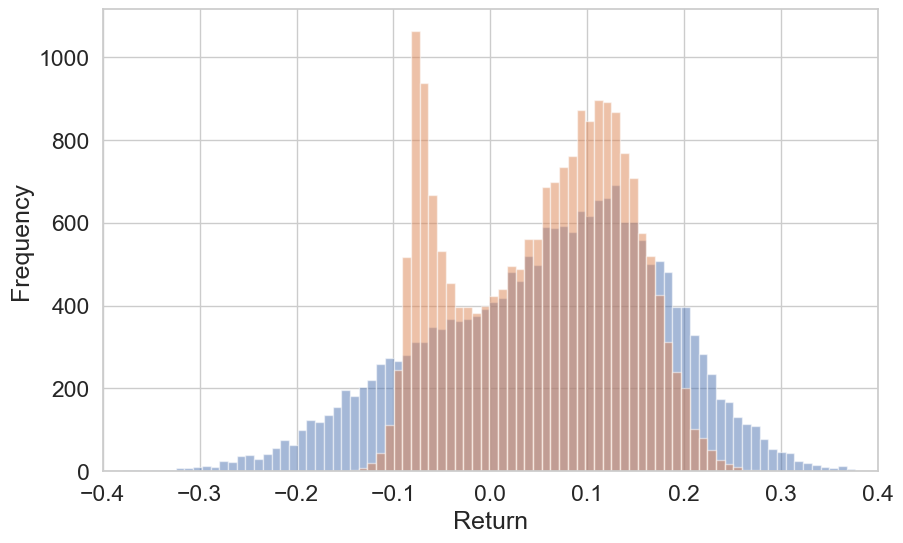}
\caption{60/40 stock/bond portfolio.}
\end{subfigure}
\caption{Histograms of simulated inflation adjusted returns of the simple
stock/bond portfolios and the collar portfolios.}
\label{fig:collar_returns}
\end{figure}

We consider the two retiree examples from \S\ref{sec:experiments}, one
upper-middle class and one lower-middle class. We use the same investment
portfolios as before, \ie, a 20/80 stock/bond portfolio in the brokerage account and a 60/40 stock/bond
portfolio in the IRA and Roth accounts, now with a self-financing collar on
the stock portfolio. We set the floor to $F=-7.5\%$. The parameter values are the
same as before, as are the realizations of the market returns, Treasury rates,
and inflation rates. A comparison between the inflation adjusted returns of the
simple stock/bond portfolios and the collar portfolios is shown in
figure~\ref{fig:collar_returns}, showing that the collar portfolios have a
narrower range of returns than the simple stock/bond portfolios, especially for
higher risk portfolios.
The annual inflation
adjusted returns of the three portfolios are found as the average over 1,000
simulations:
\[
\rho^B = 1.034, \quad \rho^I = 1.054, \quad \rho^R = 1.054.
\]

\subsection{Numerical results}
\begin{figure}
\centering
\begin{subfigure}{0.48\textwidth}
\includegraphics[width=0.9\textwidth]{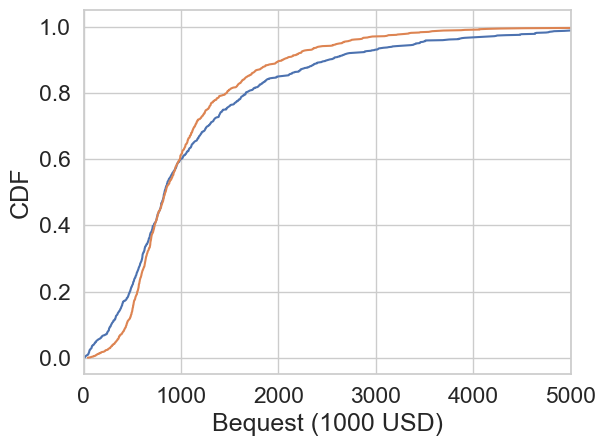}
\caption{Upper-middle-class retiree.}
\end{subfigure}
\begin{subfigure}{0.48\textwidth}
\includegraphics[width=0.9\textwidth]{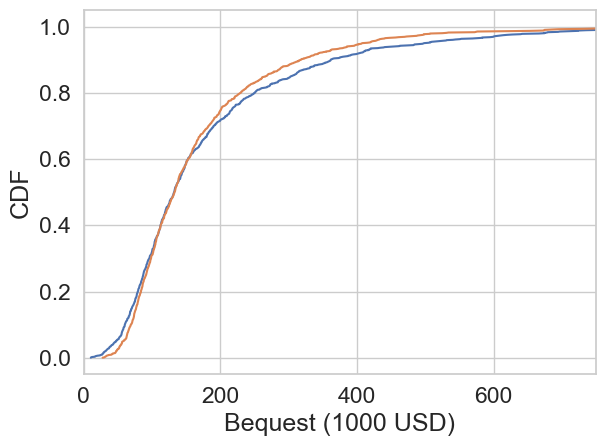}
\caption{Lower-middle-class retiree.}
\end{subfigure}
\caption{Empirical CDFs of the bequest under the MPC policy, with and without a collar, for the two retirees.}
\label{fig:cdf_bequest_collar}
\end{figure}

Figure \ref{fig:cdf_bequest_collar} shows the empirical CDSs of the bequest
under the MPC policy, with and without a collar, for the two retirees. The
collar strategy leads to a more concentrated distribution, with reduced tail
probabilities on both the left and right, compared to a simple stock/bond
portfolio. This suggests that the collar provides some protection against
extreme outcomes. Additionally, consumption remains close to the target level in
nearly all realizations when the collar is used.
The benefit of the collar is more pronounced for the upper-middle-class retiree,
where the left tail of the bequest distribution is significantly reduced. For
the lower-middle-class retiree, the effect is smaller, and it is not immediately
clear whether the collar strategy provides a meaningful improvement over a
simple stock/bond portfolio. The exact benefit depends on factors such as the
risk level of the underlying portfolio, as well as the floor of the collar. For
more aggressive allocations, say, an 80/20 stock/bond mix in the IRA and Roth,
the risk reduction from the collar is more substantial. More generally,
incorporating a collar provides a mechanism for a risk-averse retiree to
mitigate downside risk by giving up some upside potential.

\end{document}